\documentclass[11pt]{article}

\usepackage[latin1]{inputenc}\usepackage{epsfig}\usepackage{amsfonts}

\title{
  {\huge Fourier Theory on the Complex Plane III} \\
  Low-Pass Filters, Singularity Splitting \\
  and Infinite-Order Filters }

\author{
  \Large Jorge L. deLyra \\
  Department of Mathematical Physics \\
  Physics Institute \\
  University of São Paulo }

\date{March 4, 2015}

\newcommand{\ii}{\mbox{\boldmath$\imath$}}

\newcommand{\FFrac}[2]{{\displaystyle\frac{\displaystyle #1}{\displaystyle #2}}}

\newcommand{\at}[2]{\left.\rule{0em}{3ex}\right[_{\,#1}^{\,#2}}

\newcommand{\e}[1]{\,{\rm e}^{#1}}

\newcommand{\ldot}{\mbox{\Large$\cdot$}\!}

\begin{document}\maketitle

\begin{abstract}
  \noindent
  When Fourier series are employed to solve partial differential
  equations, low-pass filters can be used to regularize divergent series
  that may appear. In this paper we show that the linear low-pass filters
  defined in a previous paper can be interpreted in terms of the
  correspondence between Fourier Conjugate (FC) pairs of Definite Parity
  (DP) Fourier series and inner analytic functions, which was established
  in earlier papers. The action of the first-order linear low-pass filter
  corresponds to an operation in the complex plane that we refer to as
  ``singularity splitting'', in which any given singularity of an inner
  analytic function on the unit circle is replaced by two softer
  singularities on that same circle, thus leading to corresponding DP
  Fourier series with better convergence characteristics. Higher-order
  linear low-pass filters can be easily defined within the unit disk of
  the complex plane, in terms of the first-order one. The construction of
  infinite-order filters, which always result in $C^{\infty}$ real
  functions over the unit circle, and in corresponding DP Fourier series
  which are absolutely and uniformly convergent to these functions, is
  presented and discussed.
\end{abstract}

\section{Introduction}

In a previous paper~\cite{FTotCPI} a one-to-one correspondence between FC
(Fourier Conjugate) pairs of DP (Definite Parity) Fourier series and inner
analytic functions on the open unit disk was established. In a subsequent
paper~\cite{FTotCPII} the questions related to the convergence of such
series were examined in the light of this correspondence. In those papers
certain techniques were presented for the recovery of the real functions
from the coefficients of their DP Fourier series, which work even if the
series are divergent. This included a technique we called ``singularity
factorization'' that from the (possibly divergent) DP Fourier series of a
given real function leads to certain expressions involving alternative
trigonometric series, with better convergence characteristics, that
converge to that same real function. The reader is referred to those
papers for many of the concepts and notations used in this paper.

In another previous paper~\cite{LPFFSaPDE} certain low-pass filters acting
in the space of integrable real functions were introduced, and their use
for the regularization of divergent Fourier series in boundary value
problems was discussed. In the present paper we will show that these
low-pass filters can be interpreted and realized within the open unit disk
of the complex plane in a very simple way, in the context of the
correspondence between FC pairs of DP Fourier series and inner analytic
functions within that disk, which was discussed in the aforementioned
earlier papers~\cite{FTotCPI} and~\cite{FTotCPII}. In line with our
discussion in those papers, about ways of recovering the real functions
from their Fourier coefficients when the Fourier series do not converge,
or converge poorly, here we will show that the use of low-pass filters can
be interpreted as one more such technique. However, unlike the previous
ones it involves a certain type of approximation, and its application
changes the series and functions in a specific way, that is small in a
certain sense, as described in~\cite{LPFFSaPDE}. From a purely
mathematical standpoint the discussion of these filters consists of the
examination of the properties of a certain set of integral operators
acting in the space of integrable real functions.

Let us review briefly the facts about the filters, when defined on a
periodic interval. As given in~\cite{LPFFSaPDE}, the first-order linear
low-pass filter is defined in the following way, if we adopt as the domain
of our real functions the periodic interval $[-\pi,\pi]$. Given a real
function $f(\theta)$ of the real angular variable $\theta$ in that
interval, of which we require no more than that it be integrable, we
define from it a filtered function $f_{\epsilon}^{(1)}(\theta)$ as

\begin{equation}\label{realfilt}
  f_{\epsilon}^{(1)}(\theta)
  =
  \frac{1}{2\epsilon}
  \int_{\theta-\epsilon}^{\theta+\epsilon}d\theta'\,
  f\!\left(\theta'\right),
\end{equation}

\noindent
where the angular parameter $\epsilon\leq\pi$ is a strictly positive real
parameter which we will refer to as the {\em range} of the filter. One can
also define $f_{0}^{(1)}(\theta)$ by continuity, as the $\epsilon\to 0$
limit of this expression. The filter can be understood as a linear
integral operator acting in the space of integrable real functions, as is
done in~\cite{LPFFSaPDE}. It may be written as an integral over the whole
periodic interval involving a kernel
$K_{\epsilon}^{(1)}\!\left(\theta-\theta'\right)$ with compact support,

\begin{displaymath}
  f_{\epsilon}^{(1)}(\theta)
  =
  \int_{-\pi}^{\pi}d\theta'\,
  K_{\epsilon}^{(1)}\!\left(\theta-\theta'\right)
  f\!\left(\theta'\right),
\end{displaymath}

\noindent
where the kernel is defined as
$K_{\epsilon}^{(1)}\!\left(\theta-\theta'\right)=1/(2\epsilon)$ for
$\left|\theta-\theta'\right|<\epsilon$, and as
$K_{\epsilon}^{(1)}\!\left(\theta-\theta'\right)=0$ for
$\left|\theta-\theta'\right|>\epsilon$. Since we are in the periodic
interval, it should be noted that what we mean here by ``compact support''
is the fact that the kernel is different from zero only within an interval
contained in the periodic interval. We may have at most that the two
intervals coincide, with $\epsilon=\pi$, and in general we will assume
that we have $\epsilon\leq\pi$. The most interesting case, however, is
that in which we have $\epsilon\ll\pi$. We may say then that this kernel
is a discontinuous even function of $\left(\theta-\theta'\right)$ that has
unit integral and compact support. As shown in~\cite{LPFFSaPDE}, it can be
expressed in terms of a point-wise convergent Fourier series,

\begin{displaymath}
  K_{\epsilon}^{(1)}\!\left(\theta-\theta'\right)
  =
  \frac{1}{2\pi}
  +
  \frac{1}{\pi}
  \sum_{k=1}^{\infty}
  \left[
    \frac{\sin(k\epsilon)}{(k\epsilon)}
  \right]
  \cos\!\left[k\left(\theta-\theta'\right)\right].
\end{displaymath}

\noindent
The filter defined above has several interesting properties, which are the
reasons for its usefulness, the most important and basic ones of which are
listed and demonstrated in~\cite{LPFFSaPDE}. In this paper we will refer
to and use these properties as the occasion arises. Also, as part of the
demonstrations discussed in Section~\ref{SECinffilt} we will have the
opportunity to examine some more of these properties in
Appendix~\ref{APPInfOrd}.

As discussed in~\cite{LPFFSaPDE}, since the first-order filter defined
here is a linear operator, one can construct higher-order filters by
simply applying it multiple times to a given real function. This leads
directly to the definition of higher-order filters, for example the
second-order one, with range $2\epsilon$, and assuming that
$\epsilon\leq\pi/2$,

\begin{displaymath}
  f_{2\epsilon}^{(2)}(\theta)
  =
  \int_{-\infty}^{\infty}d\theta'\,
  K_{2\epsilon}^{(2)}\!\left(\theta-\theta'\right)
  f\!\left(\theta'\right),
\end{displaymath}

\noindent
where, as a consequence of the definition of the first-order filter, the
second-order kernel with range $2\epsilon$ is given by the application of
the first-order filter to the first-order kernel,

\begin{displaymath}
  K_{2\epsilon}^{(2)}\!\left(\theta-\theta''\right)
  =
  \int_{-\infty}^{\infty}d\theta'\,
  K_{\epsilon}^{(1)}\!\left(\theta-\theta'\right)
  K_{\epsilon}^{(1)}\!\left(\theta'-\theta''\right).
\end{displaymath}

\noindent
This second-order kernel is a continuous but non-differentiable even
function of $\left(\theta-\theta'\right)$. Due to the properties of the
first-order filter regarding its action on Fourier
expansions~\cite{LPFprop11,LPFprop12}, the second-order kernel is also
given by the absolutely and uniformly convergent Fourier series

\begin{displaymath}
  K_{2\epsilon}^{(2)}\!\left(\theta-\theta'\right)
  =
  \frac{1}{2\pi}
  +
  \frac{1}{\pi}
  \sum_{k=1}^{\infty}
  \left[
    \frac{\sin(k\epsilon)}{(k\epsilon)}
  \right]^{2}
  \cos\!\left[k\left(\theta-\theta'\right)\right],
\end{displaymath}

\noindent
so long as $\epsilon\leq\pi/2$. Both the first and second-order kernels
are even functions of $\left(\theta-\theta'\right)$ with unit integral and
compact support. The range of the first-order filter is given by
$\epsilon$, and if one just applies the filter twice as we did here, that
range doubles do $2\epsilon$. However, one may compensate for this by
simply applying twice the first-order filter with parameter $\epsilon/2$,
thus resulting in a second-order filter with range $\epsilon$, given by
the absolutely and uniformly convergent Fourier series

\begin{displaymath}
  K_{\epsilon}^{(2)}\!\left(\theta-\theta'\right)
  =
  \frac{1}{2\pi}
  +
  \frac{1}{\pi}
  \sum_{k=1}^{\infty}
  \left[
    \frac{\sin(k\epsilon/2)}{(k\epsilon/2)}
  \right]^{2}
  \cos\!\left[k\left(\theta-\theta'\right)\right],
\end{displaymath}

\noindent
so long as $\epsilon\leq\pi$. This procedure can be iterated $N$ times to
produce an order-$N$ filter with range $N\epsilon$. Given the properties
of the first-order filter regarding its action on Fourier
expansions~\cite{LPFprop11,LPFprop12}, the Fourier representation of the
order-$N$ kernel can easily be written explicitly,

\begin{equation}\label{defKordNrangeNe}
  K_{N\epsilon}^{(N)}\!\left(\theta-\theta'\right)
  =
  \frac{1}{2\pi}
  +
  \frac{1}{\pi}
  \sum_{k=1}^{\infty}
  \left[
    \frac{\sin(k\epsilon)}{(k\epsilon)}
  \right]^{N}
  \cos\!\left[k\left(\theta-\theta'\right)\right],
\end{equation}

\noindent
so long as $\epsilon\leq\pi/N$. This definition can be extended down to
the case of the order-zero kernel, with $N=0$, which is simply the Dirac
delta ``function'', and which is in fact given, as shown
in~\cite{FTotCPII}, by the divergent Fourier series

\noindent
\begin{eqnarray*}
  \delta\!\left(\theta-\theta'\right)
  & = & 
  K_{0}^{(0)}\!\left(\theta-\theta'\right)
  \\
  & = & 
  \frac{1}{2\pi}
  +
  \frac{1}{\pi}
  \sum_{k=1}^{\infty}
  \cos\!\left[k\left(\theta-\theta'\right)\right].
\end{eqnarray*}

\noindent
This can be understood as the kernel of an order-zero filter, which is the
identity almost everywhere. If we simply exchange $\epsilon$ by
$\epsilon/N$ in the expression in Equation~(\ref{defKordNrangeNe}) we get
the order-$N$ filter with range $\epsilon$, written in terms of its
Fourier expansion,

\begin{displaymath}
  K_{\epsilon}^{(N)}\!\left(\theta-\theta'\right)
  =
  \frac{1}{2\pi}
  +
  \frac{1}{\pi}
  \sum_{k=1}^{\infty}
  \left[
    \frac{\sin(k\epsilon/N)}{(k\epsilon/N)}
  \right]^{N}
  \cos\!\left[k\left(\theta-\theta'\right)\right],
\end{displaymath}

\noindent
so long as $\epsilon\leq\pi$. Note that this series converges ever faster
as $N$ increases, and that it can be differentiated $N-2$ times still
resulting in absolutely and uniformly convergent series, and $N-1$ times
still resulting in point-wise convergent series. The series for
$K_{0}^{(0)}\!\left(\theta-\theta'\right)$ is the only one which is not
convergent, and of the remaining ones that for
$K_{\epsilon}^{(1)}\!\left(\theta-\theta'\right)$ is the only one which is
not absolutely or uniformly convergent, although it is point-wise
convergent. For $N\geq 2$ all the Fourier series of the kernels,
regardless of range, are absolutely and uniformly convergent to functions
which are $C^{N-2}$ everywhere. All these kernels, regardless of order or
range, are even functions of $\left(\theta-\theta'\right)$ with unit
integral and compact support, so long as $\epsilon\leq\pi$.

Therefore, one is led to think of the possibility that in the limit
$N\to\infty$ this sequence of order-$N$ kernels with constant range
$\epsilon$ could converge to a $C^{\infty}$ kernel function
$K_{\epsilon}^{(\infty)}\!\left(\theta-\theta'\right)$ with compact
support. The corresponding infinite-order filter would then map any merely
integrable function to a $C^{\infty}$ function. Although it turns out to
be possible to construct a infinite-order kernel with such a property, it
is {\em not} to be obtained by the limit described here, as we will see
later in Section~\ref{SECinffilt}.

\section{The Low-Pass Filter on the Complex Plane}

According to the correspondence established in~\cite{FTotCPI}, to each FC
pair of DP Fourier series corresponds an inner analytic function $w(z)$
within the open unit disk. Each operation performed on the DP Fourier
series corresponds to a related operation on the inner analytic function,
possibly represented by its Taylor series around the origin. For example,
differentiation of the DP Fourier series with respect to their real
variable $\theta$ corresponds to logarithmic differentiation of $w(z)$
with respect to $z$, as shown in~\cite{FTotCPII}. If we imagine that the
first-order low-pass filter is to be implemented on the DP real functions
$f_{\rm c}(\theta)$ and $f_{\rm s}(\theta)$ associated to the DP Fourier
series, where for $z=\rho\exp(\ii\theta)$ and $\rho=1$ we have

\begin{displaymath}
  w(z)
  =
  f_{\rm c}(\theta)+\ii f_{\rm s}(\theta),
\end{displaymath}

\noindent
then it is clear that a corresponding filtering operation over $w(z)$ must
exist within the open unit disk. In this section we will give the
definition of this filtering operation on the complex plane, and derive
some of its properties.

\begin{figure}[ht]
  \centering
  \fbox{
    \epsfig{file=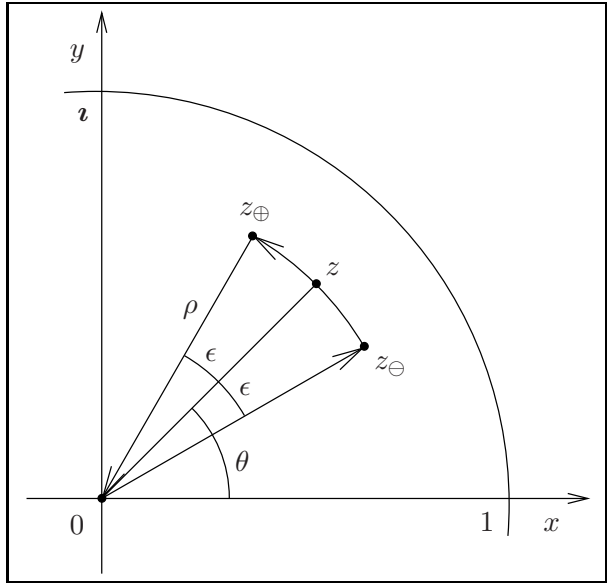,scale=1.0,angle=0}
  }
  \caption{Illustration of the definition of the first-order linear
    low-pass filter within the unit disk of the complex plane. The average
    is taken over the arc of circle from $z_{\ominus}$ to $z_{\oplus}$.}
  \label{Fig01}
\end{figure}

Consider then an inner analytic function $w(z)$, with
$z=\rho\exp(\ii\theta)$ and $0\leq\rho\leq 1$. We define from it the
corresponding filtered complex function, using the real angular range
parameter $0<\epsilon\leq\pi$, by

\begin{equation}\label{compfilt}
  w_{\epsilon}(z)
  =
  -\,
  \frac{\ii}{2\epsilon}
  \int_{z_{\ominus}}^{z_{\oplus}}dz'\,
  \frac{1}{z'}\,
  w\!\left(z'\right),
\end{equation}

\noindent
involving an integral over the arc of circle illustrated in
Figure~\ref{Fig01}, where the two extremes are given by

\noindent
\begin{eqnarray*}
  z_{\ominus}
  & = &
  z\e{-\ii\epsilon}
  \\
  & = &
  \rho\e{\ii(\theta-\epsilon)},
  \\
  z_{\oplus}
  & = &
  z\e{\ii\epsilon}
  \\
  & = &
  \rho\e{\ii(\theta+\epsilon)}.
\end{eqnarray*}

\noindent
It is important to observe that this definition can be implemented at all
the points of the unit disk, with the single additional proviso that at
$z=0$ the filter be defined as the identity. Note that the definition in
Equation~(\ref{compfilt}) has the form of a logarithmic integral, which is
the inverse operation to the logarithmic derivative, as defined and
discussed in~\cite{FTotCPII}. What we are doing here is to map the value
of the function $w(z)$ at $z$ to the average of $w(z)$ over the symmetric
arc of circle of angular length $2\epsilon$ around $z$, with constant
$\rho$. This defines a new complex function $w_{\epsilon}(z)$ at that
point. Since on the arc of circle we have that $z'=\rho\exp(\ii\theta')$
and hence that $dz'=\ii z'd\theta'$, we may also write the definition as

\begin{displaymath}
  w_{\epsilon}(z)
  =
  \frac{1}{2\epsilon}
  \int_{\theta-\epsilon}^{\theta+\epsilon}d\theta'\,
  w\!\left(\rho,\theta'\right),
\end{displaymath}

\noindent
which makes the averaging process explicitly clear. As one might expect,
just as the logarithmic differentiation of inner analytic functions
corresponds to derivatives with respect to $\theta$, the logarithmic
integration corresponds to integrals on $\theta$. Note that for
$\epsilon=\pi$ the complex filtered function $w_{\epsilon}(z)$ is simply a
constant function, possibly with removable singularities on the unit
circle. Since our real functions here, being the real or imaginary parts
of inner analytic functions, are zero-average functions, that constant is
actually zero, for all inner analytic functions.

Let us show that $w_{\epsilon}(z)$ is an inner analytic function, as
defined in~\cite{FTotCPI}. Note that since $w(z)$ is an inner analytic
function, it has the property that $w(0)=0$. Therefore we see that because
$w(0)=0$ the integrand in Equation~(\ref{compfilt}) is analytic within the
open unit circle, if defined by continuity at $z=0$. Consider therefore
the integral over the closed oriented circuit shown in Figure~\ref{Fig01},

\begin{equation}\label{compfilt2}
  \int_{0}^{z_{\ominus}}dz'\,
  \frac{1}{z'}\,
  w\!\left(z'\right)
  +
  \int_{z_{\ominus}}^{z_{\oplus}}dz'\,
  \frac{1}{z'}\,
  w\!\left(z'\right)
  +
  \int_{z_{\oplus}}^{0}dz'\,
  \frac{1}{z'}\,
  w\!\left(z'\right)
  =
  0.
\end{equation}

\noindent
Since the contour is closed and the integrand is analytic on it and within
it, this integral is zero due to the Cauchy-Goursat theorem. It follows
that we have

\begin{displaymath}
  \int_{z_{\ominus}}^{z_{\oplus}}dz'\,
  \frac{1}{z'}\,
  w\!\left(z'\right)
  =
  \int_{0}^{z_{\oplus}}dz'\,
  \frac{1}{z'}\,
  w\!\left(z'\right)
  -
  \int_{0}^{z_{\ominus}}dz'\,
  \frac{1}{z'}\,
  w\!\left(z'\right).
\end{displaymath}

\noindent
These last two integrals give the logarithmic primitive of $w(z)$ at the
two ends of the arc, as defined in~\cite{FTotCPII}. According to that
definition the logarithmic primitive of $w(z)$ is given by

\begin{displaymath}
  w^{-1\!\ldot}(z)
  =
  \int_{0}^{z}dz'\,
  \frac{1}{z'}\,
  w\!\left(z'\right),
\end{displaymath}

\noindent
where we are using the notation for the logarithmic primitive introduced
in that paper. The logarithmic primitive $w^{-1\!\ldot}(z)$ is an inner
analytic function within the open unit disk, as shown in~\cite{FTotCPII}.
It follows that we have

\begin{displaymath}
  \int_{z_{\ominus}}^{z_{\oplus}}dz'\,
  \frac{1}{z'}\,
  w\!\left(z'\right)
  =
  w^{-1\!\ldot}(z_{\oplus})
  -
  w^{-1\!\ldot}(z_{\ominus}).
\end{displaymath}

\noindent
Since the logarithmic primitive $w^{-1\!\ldot}(z)$ is an inner analytic
function, and since the functions $z_{\ominus}(z)$ and $z_{\oplus}(z)$ are
also rotated inner analytic functions, as can easily be verified, it is
reasonable to think that the right-hand side of this equation is an inner
analytic function. We have therefore for the filtered complex function

\begin{equation}\label{compfilt3}
  w_{\epsilon}(z)
  =
  -\,
  \frac{\ii}{2\epsilon}
  \left[
    w^{-1\!\ldot}(z_{\oplus})
    -
    w^{-1\!\ldot}(z_{\ominus})
  \right],
\end{equation}

\noindent
which indicates that $w_{\epsilon}(z)$ is an inner analytic function as
well. In fact, the analyticity of $w_{\epsilon}(z)$ is evident, since it
is a linear combination of two analytic functions within the open unit
disk. We must also show that $w_{\epsilon}(0)=0$ and that
$w_{\epsilon}(z)$ reduces to a real function on the interval $(-1,1)$ of
the real axis, which are the additional properties defining an inner
analytic function, as given in~\cite{FTotCPI}. It is easy to check
directly that $w_{\epsilon}(0)=0$, since $w^{-1\!\ldot}(0)=0$, given that
the logarithmic primitive is an inner analytic function. In order to
establish the remaining property, we replace $z$ by a real $x$ in the
filtered function, and taking the complex conjugate of that function with
argument $x$ we get

\begin{displaymath}
  w_{\epsilon}^{*}(x)
  =
  \frac{\ii}{2\epsilon}
  \left[
    w^{-1\!\ldot}\!\left(x\e{\ii\epsilon}\right)
    -
    w^{-1\!\ldot}\!\left(x\e{-\ii\epsilon}\right)
  \right]^{*}.
\end{displaymath}

\noindent
Now, since $w^{-1\!\ldot}(z)$ is an inner analytic function, it follows
that $w^{-1\!\ldot}(x)$ is a real function. Therefore the only relevant
participation of the number $\ii$ in the quantity within the brackets in
the expression above is that introduced explicitly via the arguments. We
have therefore, taking the complex conjugates on the right-hand side,

\noindent
\begin{eqnarray*}
  w_{\epsilon}^{*}(x)
  & = &
  \frac{\ii}{2\epsilon}
  \left[
    w^{-1\!\ldot}\!\left(x\e{-\ii\epsilon}\right)
    -
    w^{-1\!\ldot}\!\left(x\e{\ii\epsilon}\right)
  \right]
  \\
  & = &
  -\,
  \frac{\ii}{2\epsilon}
  \left[
    w^{-1\!\ldot}\!\left(x\e{\ii\epsilon}\right)
    -
    w^{-1\!\ldot}\!\left(x\e{-\ii\epsilon}\right)
  \right]
  \\
  & = &
  w_{\epsilon}(x),
\end{eqnarray*}

\noindent
so that $w_{\epsilon}(z)$ reduces to a real function on the interval
$(-1,1)$ of the real axis. This establishes, therefore, that the filtered
complex function $w_{\epsilon}(z)$ is in fact an inner analytic function.
In addition to this, since logarithmic integration softens the
singularities of $w(z)$ by one degree, as discussed in~\cite{FTotCPII}, we
see that $w_{\epsilon}(z)$ will have all its singularities softened by one
degree as compared to those of $w(z)$.

Observe that, if we take the limit $\rho\to 1$ to the unit circle in such
a way that $z$ tends to a singularity of $w(z)$ at the position $\theta$,
it immediately follows that $w_{\epsilon}(z)$ has two singularities, each
softened by one degree, at the positions $(\theta-\epsilon)$ and
$(\theta+\epsilon)$. What we have here is what we will refer to as the
process of {\em singularity splitting}, for we see that the application of
the filter has the effect of interchanging a harder singularity at
$\theta$ by two softer singularities at $(\theta-\epsilon)$ and
$(\theta+\epsilon)$. In particular, this will always decrease the degree
of hardness, or increase the degree of softness, of all the dominant
singularities on the unit circle, by one degree. This in turn is important
because the dominant singularities determine the level and mode of
convergence of the DP Fourier series, as discussed in~\cite{FTotCPII}.

Observe that the filtering operation does {\em not} stay within a single
integral-differential chain of inner analytic functions, as defined
in~\cite{FTotCPII}, since it changes the location of the singularities of
the inner analytic function it is applied to. Instead, it passes to
another such chain, while at the same time changing to the next link in
the new chain, in the softening direction, since it softens the
singularities by one degree. The new function reached in this way is not
directly related to the original one by straight logarithmic
integration. The new function is, however, close to the original function,
so long as $\epsilon$ is small, according to a criterion that has a clear
physical meaning, as explained in~\cite{LPFFSaPDE}.

Since the complex-plane definition of the first-order low-pass filter in
the open unit disk reproduces the definition of the filter as given in
Equation~(\ref{realfilt}) directly in terms of the corresponding real
functions on the unit circle, it also reproduces all the properties of the
filter when acting on the real functions, which were discussed and
demonstrated in~\cite{LPFFSaPDE}. In some cases there are corresponding
properties of the filter in terms of the complex functions. By
construction it is clear that, just as $w(z)$, the function
$w_{\epsilon}(z)$ is periodic in $\theta$, with period $2\pi$, which is a
generalization to the complex plane of one of the properties of the
filter~\cite{LPFprop10}. In addition to this, since it acts on inner
analytic functions, which are analytic and thus always continuous and
differentiable, it is quite clear that the filter becomes the identity
operation in the $\epsilon\to 0$ limit. We can see this from the
complex-plane definition in Equation~(\ref{compfilt3}). If we consider the
variation of $\theta$ between $z_{\oplus}$ and $z_{\ominus}$, which is
given in terms of the parameter $\epsilon$ by $\delta\theta=2\epsilon$,
and we take the $\epsilon\to 0$ limit of that expression, we get

\noindent
\begin{eqnarray*}
  \lim_{\epsilon\to 0}
  w_{\epsilon}(z)
  & = &
  -\ii\,
  \lim_{\epsilon\to 0}
  \frac{w^{-1\!\ldot}(z_{\oplus})-w^{-1\!\ldot}(z_{\ominus})}{2\epsilon}
  \\
  & = &
  -\ii\,
  \lim_{\delta\theta\to 0}
  \frac{w^{-1\!\ldot}(z_{\oplus})-w^{-1\!\ldot}(z_{\ominus})}{\delta\theta}
  \\
  & = &
  z
  \lim_{\delta z\to 0}
  \frac{w^{-1\!\ldot}(z_{\oplus})-w^{-1\!\ldot}(z_{\ominus})}{\delta z},
\end{eqnarray*}

\noindent
where we used the fact that in the limit $\delta z=\ii z\delta\theta$.
The limit above defines the logarithmic derivative, so that we have

\noindent
\begin{eqnarray*}
  \lim_{\epsilon\to 0}
  w_{\epsilon}(z)
  & = &
  z\,
  \frac{d}{dz}w^{-1\!\ldot}(z)
  \\
  & = &
  w(z),
\end{eqnarray*}

\noindent
since we have the logarithmic derivative of the logarithmic primitive, and
the operations of logarithmic differentiation and logarithmic integration
are the inverses of one another. This establishes that in the $\epsilon\to
0$ limit the filter becomes the identity when acting on the inner analytic
functions, which is a generalization to the complex plane of another
property of the filter~\cite{LPFprop06}. In fact, this property within the
open unit disk is somewhat stronger than the corresponding property on the
unit circle, since in this case we have exactly the identity in all cases,
while in the real case we had only the identity almost everywhere.

Taken in the light of the classification of singularities and modes of
convergence which was given in~\cite{FTotCPII}, we can see immediately the
consequences of this process of singularity splitting on the mode of
convergence of the DP Fourier series associated to the inner analytic
function, and on the analytic character of the corresponding DP real
functions. Let us recall from the earlier papers~\cite{FTotCPI}
and~\cite{FTotCPII} that given an inner analytic function

\begin{displaymath}
  w(z)
  =
  f_{\rm c}(\rho,\theta)+\ii f_{\rm s}(\rho,\theta),
\end{displaymath}

\noindent
where $z=\rho\exp(\ii\theta)$, and its Taylor series around $z=0$,

\begin{displaymath}
  w(z)
  =
  \sum_{k=1}^{\infty}
  a_{k}z^{k},
\end{displaymath}

\noindent
which is convergent at least on the open unit disk, it follows that on the
unit circle we have the real functions $f_{\rm c}(\theta)=f_{\rm
  c}(1,\theta)$ and $f_{\rm s}(\theta)=f_{\rm s}(1,\theta)$, associated to
the FC pair of DP Fourier series

\noindent
\begin{eqnarray*}
  f_{\rm c}(\theta)
  & = &
  \sum_{k=1}^{\infty}
  a_{k}\cos(k\theta),
  \\
  f_{\rm s}(\theta)
  & = &
  \sum_{k=1}^{\infty}
  a_{k}\sin(k\theta).
\end{eqnarray*}

\noindent
After the action of the filter we have corresponding relations for the
filtered functions,

\noindent
\begin{eqnarray*}
  w_{\epsilon}(z)
  & = &
  f_{\epsilon,{\rm c}}(\rho,\theta)+\ii f_{\epsilon,{\rm s}}(\rho,\theta)
  \\
  & = &
  \sum_{k=1}^{\infty}
  a_{\epsilon,k}z^{k},
  \\
  f_{\epsilon,{\rm c}}(\theta)
  & = &
  \sum_{k=1}^{\infty}
  a_{\epsilon,k}\cos(k\theta),
  \\
  f_{\epsilon,{\rm s}}(\theta)
  & = &
  \sum_{k=1}^{\infty}
  a_{\epsilon,k}\sin(k\theta).
\end{eqnarray*}

\noindent
The results obtained in~\cite{FTotCPII} relate the nature of the dominant
singularities of $w(z)$ on the unit circle with the mode of convergence of
the corresponding DP Fourier series and with the analytical character of
the corresponding DP real functions $f_{\rm c}(\theta)$ and $f_{\rm
  s}(\theta)$, for a large class of inner analytic functions and
corresponding DP real functions. The same relations also hold for
$w_{\epsilon}(z)$, of course. Assuming that the functions at issue here
are within that class, we may derive some of the properties of the
first-order low-pass filter, as defined on the complex plane.

For one thing, if the original real functions are continuous, then
according to the classification introduced in~\cite{FTotCPII} the original
inner analytic function has dominant singularities that are soft, with any
degree of softness starting from borderline soft singularities (that is,
with degree of softness zero), and the DP Fourier series are absolutely
and uniformly convergent. In this case the action of the filter results in
an inner analytic function with dominant singularities that have a degree
of softness equal to $1$ or larger, thus implying that the corresponding
filtered real functions are differentiable. We thus reproduce in the
complex formalism one of the properties of the first-order
filter~\cite{LPFprop03}, namely that if a real function is continuous then
the corresponding filtered function is differentiable.

In addition to this, if the original real functions are integrable but not
continuous, then according to the classification introduced
in~\cite{FTotCPII} the original inner analytic function has dominant
singularities that are borderline hard ones (that is, with degree of
hardness zero), and the DP Fourier series are convergent almost
everywhere, but not absolutely or uniformly convergent. In this case the
action of the filter results in an inner analytic function with dominant
singularities which are borderline soft, thus implying that the
corresponding filtered real functions are continuous. Also, in this case
the filtered DP Fourier series become absolutely and uniformly convergent.
We thus reproduce in the complex formalism another one of the properties
of the first-order filter~\cite{LPFprop04}, namely that if a real function
is discontinuous then the corresponding filtered function is continuous.

Since the filter acts only on the variable $\theta$, some of the
properties of the filter defined on the real line, and hence on the unit
circle, are translated transparently to the complex formalism. For
example, the action on the filter on the Fourier expansions encoded into
the angular part of the complex Taylor expansions is determined by its
action on the elements of the Fourier basis, as shown
in~\cite{LPFprop11,LPFprop12}. If we apply the filter as defined in
Equation~(\ref{compfilt2}) to the functions of the basis we get

\noindent
\begin{eqnarray*}
  \frac{1}{2\epsilon}
  \int_{\theta-\epsilon}^{\theta+\epsilon}d\theta'\,
  \cos\!\left(k\theta'\right)
  & = &
  \left[
    \frac{\sin(k\epsilon)}{(k\epsilon)}
  \right]
  \cos(k\theta),
  \\
  \frac{1}{2\epsilon}
  \int_{\theta-\epsilon}^{\theta+\epsilon}d\theta'\,
  \sin\!\left(k\theta'\right)
  & = &
  \left[
    \frac{\sin(k\epsilon)}{(k\epsilon)}
  \right]
  \sin(k\theta).
\end{eqnarray*}

\noindent
This means that the elements of that basis are eigenfunctions of the
filter, interpreted as an operator. It also determines the eigenvalues,
given by the ratio shown within brackets, which is known as the sinc
function of the variable $(k\epsilon)$. What this means is that the filter
acts of an extremely simple way on the Fourier expansions. It then follows
that the same is true, of course, for the Taylor series of the
corresponding inner analytic functions. If we write the Taylor expansion
of a given inner analytic function in polar coordinates, with
$z=\rho\exp(\ii\theta)$, we get

\begin{displaymath}
  w(z)
  =
  \sum_{k=1}^{\infty}
  a_{k}
  \rho^{k}
  \left[
    \cos(k\theta)
    +
    \ii
    \sin(k\theta)
  \right],
\end{displaymath}

\noindent
and from this follows at once the corresponding expansion for the filtered
function

\begin{displaymath}
  w_{\epsilon}(z)
  =
  \sum_{k=1}^{\infty}
  a_{k}
  \left[
    \frac{\sin(k\epsilon)}{(k\epsilon)}
  \right]
  \rho^{k}
  \left[
    \cos(k\theta)
    +
    \ii
    \sin(k\theta)
  \right].
\end{displaymath}

\noindent
What this means is that the Taylor coefficients $a_{\epsilon,k}$ of
$w_{\epsilon}(z)$ are given by

\begin{displaymath}
  a_{\epsilon,k}
  =
  \left[
    \frac{\sin(k\epsilon)}{(k\epsilon)}
  \right]
  a_{k},
\end{displaymath}

\noindent
in terms of the Taylor coefficients $a_{k}$ of $w(z)$, a fact that can be
shown directly from the definition of the coefficients, as one can see in
Section~\ref{APPFourCoef} of Appendix~\ref{APPProofs}. If we take the
$\rho\to 1$ limit this corresponds to the filtered real functions

\noindent
\begin{eqnarray*}
  f_{\epsilon,{\rm c}}(\theta)
  & = &
  \sum_{k=1}^{\infty}
  a_{k}
  \left[
    \frac{\sin(k\epsilon)}{(k\epsilon)}
  \right]
  \cos(k\theta),
  \\
  f_{\epsilon,{\rm s}}(\theta)
  & = &
  \sum_{k=1}^{\infty}
  a_{k}
  \left[
    \frac{\sin(k\epsilon)}{(k\epsilon)}
  \right]
  \sin(k\theta).
\end{eqnarray*}

\noindent
It follows therefore that the same relation holds for the Fourier
coefficients of $f_{\epsilon,{\rm c}}(\theta)$ and $f_{\epsilon,{\rm
    s}}(\theta)$, in terms of the Fourier coefficients of $f_{\rm
  c}(\theta)$ and $f_{\rm s}(\theta)$.

\begin{figure}[ht]
  \centering
  \fbox{
    \epsfig{file=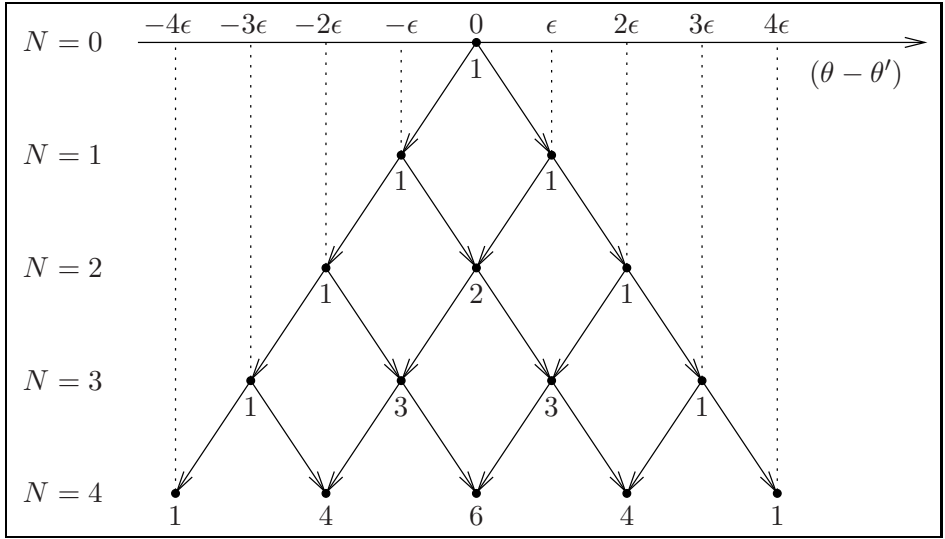,scale=1.0,angle=0}
  }
  \caption{The iteration of the first-order filter to produce an order-$N$
    filter, showing the structure of the Pascal triangle and the linear
    increase of the range. The original singularity is at $\theta'$. The
    numbers near the vertices of the triangles represent the number of
    softened singularities superposed at that point.}
  \label{Fig02}
\end{figure}

\subsection{Higher-Order Filters}

If one simply iterates $N$ times the procedure in
Equation~(\ref{compfilt3}), which is equivalent to the definition of the
first-order linear low-pass filter within the unit disk of the complex
plane, one gets the corresponding higher-order filters in the
complex-plane representation. For example, a second-order filter with
range $2\epsilon$ can be obtained by applying the first-order filter
twice, and results in the complex-plane definition

\begin{displaymath}
  w_{2\epsilon}^{(2)}(z)
  =
  \left(\frac{-\ii}{2\epsilon}\right)^{2}
  \left[
    w^{-2\ldot}\!\left(z\e{2\ii\epsilon}\right)
    -
    2w^{-2\ldot}\!\left(z\right)
    +
    w^{-2\ldot}\!\left(z\e{-2\ii\epsilon}\right)
  \right].
\end{displaymath}

\noindent
Note that only the second logarithmic primitive of $w(z)$ appears here,
and hence that all singularities are softened by two degrees. A
corresponding second-order filter with range $\epsilon$ can then be
obtained by the exchange of $\epsilon$ by $\epsilon/2$,

\begin{displaymath}
  w_{\epsilon}^{(2)}(z)
  =
  \left(\frac{-\ii}{\epsilon}\right)^{2}
  \left[
    w^{-2\ldot}\!\left(z\e{\ii\epsilon}\right)
    -
    2w^{-2\ldot}\!\left(z\right)
    +
    w^{-2\ldot}\!\left(z\e{-\ii\epsilon}\right)
  \right].
\end{displaymath}

\noindent
It is now possible to define an $N^{\rm th}$ order filter, with range
$N\epsilon$, by iterating this procedure repeatedly. If we look at it in
terms of the singularities on the unit circle, the iteration corresponds
to recursive singularity splitting as shown in Figure~\ref{Fig02}. In this
diagram we can see the structure of the Pascal triangle, and the linear
increase of the resulting range with $N$. Note that at the $N^{\rm th}$
iteration there are $N+1$ softened singularities within the interval
$[-N\epsilon,N\epsilon]$ of the variable $\left(\theta-\theta'\right)$. It
is important to observe that, while the singularities become progressively
softer as one goes down the diagram, it is also the case that more and
more singularities are superposed at the same point, particularly near the
central vertical line of the triangle. From the structure of the Pascal
triangle the coefficients of the superposition are easily obtained, so
that from this diagram it is not too difficult to obtain the expression
for the $N^{\rm th}$ order filter, with range $N\epsilon$, which turns out
to be

\begin{displaymath}
  w_{N\epsilon}^{(N)}(z)
  =
  \left(\frac{-\ii}{2\epsilon}\right)^{N}
  \sum_{n=0}^{N}
  \frac{(-1)^{n}N!}{n!(N-n)!}\,
  w^{-N\ldot}\!\left(z\e{\ii[N-2n]\epsilon}\right).
\end{displaymath}

\noindent
In this case the range of the changes introduced in the real functions by
the order-$N$ filter has the value $N\epsilon$. This means that, given a
fixed value of $\epsilon$, the iteration of the first-order filter cannot
be done indefinitely inside the periodic interval $[-\pi,\pi]$ without the
range eventually becoming larger than the period. However, one may reduce
the resulting range back to $\epsilon$ by simply using for the
construction the linear filter with range $\epsilon/N$, resulting in

\begin{displaymath}
  w_{\epsilon}^{(N)}(z)
  =
  \left(\frac{-\ii N}{2\epsilon}\right)^{N}
  \sum_{n=0}^{N}
  \frac{(-1)^{n}N!}{n!(N-n)!}\,
  w^{-N\ldot}\!\left(z\e{\ii[1-2n/N]\epsilon}\right).
\end{displaymath}

\noindent
With this renormalization of the parameter $\epsilon$ it is now possible
to do the iteration of the first-order filter indefinitely inside the
periodic interval $[-\pi,\pi]$, keeping the range constant, and therefore
to define filters of arbitrarily high orders. In this case a singularity
at $\theta'$ on the unit circle will be split into $N+1$ singularities
softened by $N$ degrees, homogeneously distributed within the interval
$[-\epsilon,\epsilon]$ of the variable $\left(\theta-\theta'\right)$. One
may even consider iterating the filter an infinite number of times in this
way, keeping the range constant. However, this does {\em not} work quite
as one might expect at first. A detailed discussion of this case can be
found in Section~\ref{SECinffilt}.

Note that since these higher-order filters are obtained by the repeated
application of the first-order one, they inherit from it many of its
properties. For example, they are all the identity when applied to linear
real functions on the unit circle~\cite{LPFprop01}, and they all maintain
the periodicity of periodic functions~\cite{LPFprop10}. Also, they all
have the elements of the Fourier basis as eigenfunctions and hence they
all commute with the second-derivative operator, as demonstrated
in~\cite{LPFFSaPDE}. In terms of the DP Fourier series, if one considers
the $N$-fold repeated application of the first-order filter to the
original real function, since each instance of the first-order filter
contributes the same factor to the coefficients, as shown
in~\cite{LPFFSaPDE}, one simply gets for the filtered real functions

\noindent
\begin{eqnarray*}
  f_{N\epsilon,{\rm c}}^{(N)}(\theta)
  & = &
  \sum_{k=1}^{\infty}
  a_{k}
  \left[
    \frac{\sin(k\epsilon)}{(k\epsilon)}
  \right]^{N}
  \cos(k\theta),
  \\
  f_{N\epsilon,{\rm s}}^{(N)}(\theta)
  & = &
  \sum_{k=1}^{\infty}
  a_{k}
  \left[
    \frac{\sin(k\epsilon)}{(k\epsilon)}
  \right]^{N}
  \sin(k\theta).
\end{eqnarray*}

\noindent
This will of course imply that the filtered DP Fourier series converge
significantly faster than the original ones, and to significantly smoother
functions. In this case the range of the changes introduced in the real
functions has the value $N\epsilon$. Once more one may reduce the
resulting range back to $\epsilon$, using the linear filter with range
$\epsilon/N$, thus leading to

\noindent
\begin{eqnarray*}
  f_{\epsilon,{\rm c}}^{(N)}(\theta)
  & = &
  \sum_{k=1}^{\infty}
  a_{k}
  \left[
    \frac{\sin(k\epsilon/N)}{(k\epsilon/N)}
  \right]^{N}
  \cos(k\theta),
  \\
  f_{\epsilon,{\rm s}}^{(N)}(\theta)
  & = &
  \sum_{k=1}^{\infty}
  a_{k}
  \left[
    \frac{\sin(k\epsilon/N)}{(k\epsilon/N)}
  \right]^{N}
  \sin(k\theta).
\end{eqnarray*}

\noindent
This modification changes only the range of the alterations introduced in
the real functions by the order-$N$ filter, and not the level of
smoothness of the resulting filtered functions, which depends only on $N$.

Since they are themselves real functions defined on a circle of radius
$\rho\leq 1$ in the complex plane, centered at the origin, the kernels of
the order-$N$ filters can also be represented by inner analytic functions
within the corresponding open disk. This is a simple extension of the
structure we developed in the earlier papers~\cite{FTotCPI}
and~\cite{FTotCPII}. If $z_{1}=\rho_{1}\exp(\ii\theta_{1})$ is a point on
the circle $\rho_{1}\leq 1$ and $z=\rho\exp(\ii\theta)$ a point inside the
corresponding disk, the kernels of constant range $\epsilon$ can be
written as the real parts of the complex kernels

\begin{displaymath}
  \kappa_{\epsilon}^{(N)}\!\left(z,z_{1}\right)
  =
  \frac{1}{2\pi}
  +
  \frac{1}{\pi}
  \sum_{k=1}^{\infty}
  \left[
    \frac{\sin(k\epsilon/N)}{(k\epsilon/N)}
  \right]^{N}
  \left(
    \frac{z}{z_{1}}
  \right)^{k},
\end{displaymath}

\noindent
where it should be noted that the coefficients are real. Except for the
constant term this is the Taylor series of an inner analytic function
inside the disk of radius $\rho_{1}$, rotated by the angle
$\theta_{1}$. If we take the limit $\rho\to\rho_{1}$ we get

\noindent
\begin{eqnarray*}
  \kappa_{\epsilon}^{(N)}(\theta-\theta_{1})
  & = &
  \frac{1}{2\pi}
  +
  \frac{1}{\pi}
  \sum_{k=1}^{\infty}
  \left[
    \frac{\sin(k\epsilon/N)}{(k\epsilon/N)}
  \right]^{N}
  \e{\ii k(\theta-\theta_{1})}
  \\
  & = &
  \frac{1}{2\pi}
  +
  \frac{1}{\pi}
  \sum_{k=1}^{\infty}
  \left[
    \frac{\sin(k\epsilon/N)}{(k\epsilon/N)}
  \right]^{N}
  \cos[k(\theta-\theta_{1})]
  +
  \\
  &   &
  \hspace{1.6em}
  +
  \ii\,
  \frac{1}{\pi}
  \sum_{k=1}^{\infty}
  \left[
    \frac{\sin(k\epsilon/N)}{(k\epsilon/N)}
  \right]^{N}
  \sin[k(\theta-\theta_{1})],
\end{eqnarray*}

\noindent
and therefore we have

\begin{displaymath}
  \Re\!\left[\kappa_{\epsilon}^{(N)}(\theta-\theta_{1})\right]
  =
  K_{\epsilon}^{(N)}(\theta-\theta_{1}).
\end{displaymath}

\noindent
Note that using this complex-plane representation it is easy to prove that
the kernels of the order-$N$ filters have unit integral. We consider the
integral over the circle $C_{1}$ of radius $\rho_{1}$, that appears in the
Cauchy integral formula for $\kappa_{\epsilon}^{(N)}(z,z_{1})$ around
$z=0$,

\begin{displaymath}
  \frac{1}{2\pi\ii}
  \oint_{C_{1}}dz\,
  \frac{1}{z}\,
  \kappa_{\epsilon}^{(N)}(z,z_{1})
  =
  \kappa_{\epsilon}^{(N)}(0,z_{1}).
\end{displaymath}

\noindent
Since we have the value $\kappa_{\epsilon}^{(N)}(0,z_{1})=1/(2\pi)$, we
get

\begin{displaymath}
  \frac{1}{\ii}
  \oint_{C_{1}}dz\,
  \frac{1}{z}\,
  \kappa_{\epsilon}^{(N)}(z,z_{1})
  =
  1.
\end{displaymath}

\noindent
If we now write the integral explicitly over the circle, with
$z=\rho_{1}\exp(\ii\theta)$ and $dz=\ii zd\theta$, we get

\begin{displaymath}
  \int_{-\pi}^{\pi}d\theta\,
  \kappa_{\epsilon}^{(N)}(\theta-\theta_{1})
  =
  1.
\end{displaymath}

\noindent
Finally, if we consider explicitly the real and imaginary parts we get

\noindent
\begin{eqnarray*}
  \int_{-\pi}^{\pi}d\theta\,
  \left\{
    \Re\!\left[\kappa_{\epsilon}^{(N)}(\theta-\theta_{1})\right]
    +
    \ii
    \Im\!\left[\kappa_{\epsilon}^{(N)}(\theta-\theta_{1})\right]
  \right\}
  & = &
  1
  \;\;\;\Rightarrow
  \\
  \int_{-\pi}^{\pi}d\theta\,
  \Re\!\left[\kappa_{\epsilon}^{(N)}(\theta-\theta_{1})\right]
  & = &
  1,
  \\
  \int_{-\pi}^{\pi}d\theta\,
  \Im\!\left[\kappa_{\epsilon}^{(N)}(\theta-\theta_{1})\right]
  & = &
  0.
\end{eqnarray*}

\noindent
Since the real part is $K_{\epsilon}^{(N)}(\theta-\theta_{1})$, the result
follows,

\begin{displaymath}
  \int_{-\pi}^{\pi}d\theta\,
  K_{\epsilon}^{(N)}(\theta-\theta_{1})
  =
  1,
\end{displaymath}

\noindent
for all $N$.

\section{The Infinite-Order Filter}\label{SECinffilt}

Let us now discuss the possibility of constructing infinite-order filters
with compact support. As was mentioned before, it would be an interesting
thing to have the definition of an infinite-order linear low-pass filter.
If we consider the linear low-pass filter of order $N$ and a fixed range
$\epsilon$, that can be described in terms of the inner-analytic functions
as

\begin{equation}\label{filterordNrangee}
  w_{\epsilon}^{(N)}(z)
  =
  \left(\frac{-\ii N}{2\epsilon}\right)^{N}
  \sum_{n=0}^{N}
  \frac{(-1)^{n}N!}{n!(N-n)!}\,
  w^{-N\ldot}\!\left[z\e{\ii(1-2n/N)\epsilon}\right],
\end{equation}

\noindent
it is natural to ask that happens if we take the $N\to\infty$ limit. This
cannot be described simply as an infinite iteration of the first-order
linear filter, since the limiting process changes the range of that filter
to zero. On the other hand, all the filtered complex functions
$w_{\epsilon}^{(N)}(z)$ exist and are inner analytic, as a sequence
indexed by $N$, for all $N$, so that it is reasonable to think that the
limit should also exist and should also be an inner analytic function, at
least inside the open unit disk of the complex plane. However, it is
important to keep in mind that it is far less clear what happens when one
takes the limit from the open unit disk to the unit circle, {\em after}
one first takes the $N\to\infty$ limit.

In this section we will endeavor to construct an infinite-order filter
with compact support. If this endeavor succeeds, then there is an
interesting consequence of the eventual construction of such an
infinite-order filter, regarding the construction of $C^{\infty}$
functions with compact support. If it turns out to be possible to define
this infinite-order filter with a finite range $\epsilon$ in terms of an
integral involving a well-defined infinite-order kernel with compact
support,

\begin{displaymath}
  f_{\epsilon}^{(\infty)}(\theta)
  =
  \int_{-\pi}^{\pi}d\theta'\,
  K_{\epsilon}^{(\infty)}\!\left(\theta-\theta'\right)
  f(\theta'),
\end{displaymath}

\noindent
then it would in principle be possible to use this filter operator to
transform any integrable function into a $C^{\infty}$ function, making
changes only within a finite range $\epsilon$ that can be as small as one
wishes. In our current case here this infinite-order kernel would be
written as the limit

\noindent
\begin{eqnarray*}
  K_{\epsilon}^{(\infty)}\!\left(\theta-\theta'\right)
  & = &
  \lim_{N\to\infty}
  K_{\epsilon}^{(N)}\!\left(\theta-\theta'\right)
  \\
  & = &
  \frac{1}{2\pi}
  +
  \frac{1}{\pi}
  \lim_{N\to\infty}
  \sum_{k=1}^{\infty}
  \left[
    \frac{\sin(k\epsilon/N)}{(k\epsilon/N)}
  \right]^{N}
  \cos\!\left[k\left(\theta-\theta'\right)\right].
\end{eqnarray*}

\noindent
The idea here is that the kernel
$K_{\epsilon}^{(\infty)}\!\left(\theta-\theta'\right)$ would then be
itself a $C^{\infty}$ function, and that due to the properties of the
first-order filter, it would also have unit integral. However, the fact is
that the limit above does not behave as one might expect at first. If we
consider the $N\to\infty$ limit of the coefficients, we have

\noindent
\begin{eqnarray*}
  \lefteqn
  {
    \lim_{N\to\infty}
    \left[
      \left(\frac{N}{k\epsilon}\right)
      \sin\!\left(\frac{k\epsilon}{N}\right)
    \right]^{N}
  }
  &   &
  \\
  & = &
  \lim_{N\to\infty}
  \left[
    \left(\frac{N}{k\epsilon}\right)
    \sum_{j=0}^{\infty}
    \frac{(-1)^{j}}{(2j+1)!}
    \left(\frac{k\epsilon}{N}\right)^{2j+1}
  \right]^{N}
  \\
  & = &
  \lim_{N\to\infty}
  \left[
    \sum_{j=0}^{\infty}
    \frac{(-1)^{j}}{(2j+1)!}
    \left(\frac{k\epsilon}{N}\right)^{2j}
  \right]^{N}
  \\
  & = &
  \lim_{N\to\infty}
  \left[
    1
    -
    \frac{1}{6}
    \left(\frac{k\epsilon}{N}\right)^{2}
    +
    \frac{1}{120}
    \left(\frac{k\epsilon}{N}\right)^{4}
    -
    \frac{1}{720}
    \left(\frac{k\epsilon}{N}\right)^{6}
    +
    \ldots
  \right]^{N}.
\end{eqnarray*}

\noindent
If one expands the power $N$, there is one term equal to $1$ and all other
terms have powers of $N$ in the denominator. If we write the terms that
have up to four powers of $N$ in the denominator, we get

\noindent
\begin{eqnarray*}
  \lefteqn
  {
    \lim_{N\to\infty}
    \left[
      \left(\frac{N}{k\epsilon}\right)
      \sin\!\left(\frac{k\epsilon}{N}\right)
    \right]^{N}
  }
  &   &
  \\
  & = &
  \lim_{N\to\infty}
  \left[
    1
    -
    \frac{1}{6}\,
    N
    \left(\frac{k\epsilon}{N}\right)^{2}
    +
    \frac{1}{36}\,
    \frac{N(N-1)}{2}
    \left(\frac{k\epsilon}{N}\right)^{4}
  \right.
  +
  \\
  &   &
  \hspace{3em}
  +
  \left.
    \frac{1}{120}\,
    N
    \left(\frac{k\epsilon}{N}\right)^{4}
    -
    \frac{1}{720}\,
    N(N-1)
    \left(\frac{k\epsilon}{N}\right)^{6}
    +
    \ldots
  \right].
\end{eqnarray*}

\noindent
A more rigorous analysis of this limit would require more careful
consideration of the convergence of this series, since in principle one
must be careful with the interchange of the $N\to\infty$ limit and the
$j\to\infty$ limit of the series. However, it turns out that this rough
discussion suffices for our purposes here. As one can see, all terms
except the first have at least one factor of $N$ in the denominator, and
therefore we should expect that

\begin{displaymath}
    \lim_{N\to\infty}
    \left[
      \left(\frac{N}{k\epsilon}\right)
      \sin\!\left(\frac{k\epsilon}{N}\right)
    \right]^{N}
    =
    1,
\end{displaymath}

\noindent
for all $k$. This implies that we have for the finite-range kernel, in the
$N\to\infty$ limit,

\noindent
\begin{eqnarray*}
  K_{\epsilon}^{(\infty)}\!\left(\theta-\theta'\right)
  & = &
  \lim_{N\to\infty}
  K_{\epsilon}^{(N)}\!\left(\theta-\theta'\right)
  \\
  & = &
  \frac{1}{2\pi}
  +
  \frac{1}{\pi}
  \sum_{k=1}^{\infty}
  \cos\!\left[k\left(\theta-\theta'\right)\right],
\end{eqnarray*}

\noindent
which is in fact the Fourier expansion of the Dirac delta ``function'',
which is something of an unexpected outcome! In other words, the limit of
this sequence of progressively smoother functions is not even a function,
but a singular object instead. This is actually very similar to the
representation of the delta ``function'' by an infinite sequence of
normalized Gaussian functions.

\begin{figure}[ht]
  \centering
  \fbox{
    \epsfig{file=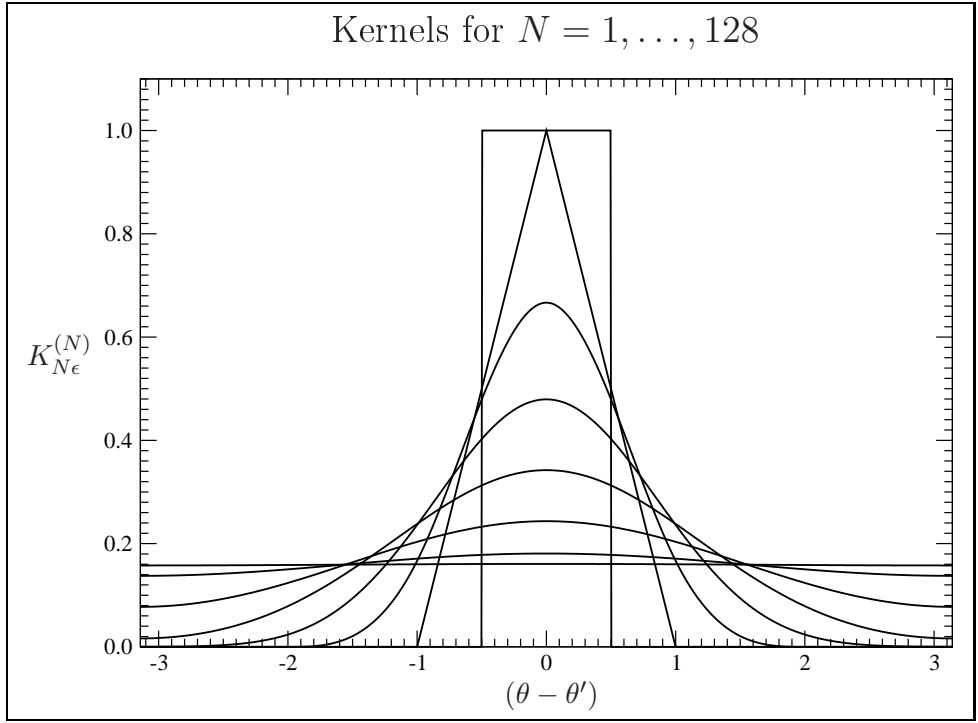,scale=1.0,angle=0}
  }
  \caption{The kernels of several filters with increasing range
    $N\epsilon$, obtained via the use of their Fourier series, for
    $\epsilon=0.5$, for values of $N$ increasing exponentially, in the set
    $\{1,2,4,8,16,32,64,128\}$, plotted as functions of
    $\left(\theta-\theta'\right)$ within the periodic interval
    $[-\pi,\pi]$.}
  \label{Fig03}
\end{figure}

A little numerical exploration is useful at this point to establish some
simple mathematical facts about these infinite-order kernels. For
completeness, let us go momentarily back to the straightforward multiple
superposition of the first-order filter. If we take the $N\to\infty$ limit
of the kernel of order $N$ with range $N\epsilon$, we might try to define
a first infinite-order kernel with infinite range as

\begin{displaymath}
  K_{\infty}^{(\infty)}\!\left(\theta-\theta'\right)
  =
  \frac{1}{2\pi}
  +
  \frac{1}{\pi}
  \lim_{N\to\infty}
  \sum_{k=1}^{\infty}
  \left[
    \frac{\sin(k\epsilon)}{(k\epsilon)}
  \right]^{N}
  \cos\!\left[k\left(\theta-\theta'\right)\right].
\end{displaymath}

\begin{figure}[ht]
  \centering
  \fbox{
    \epsfig{file=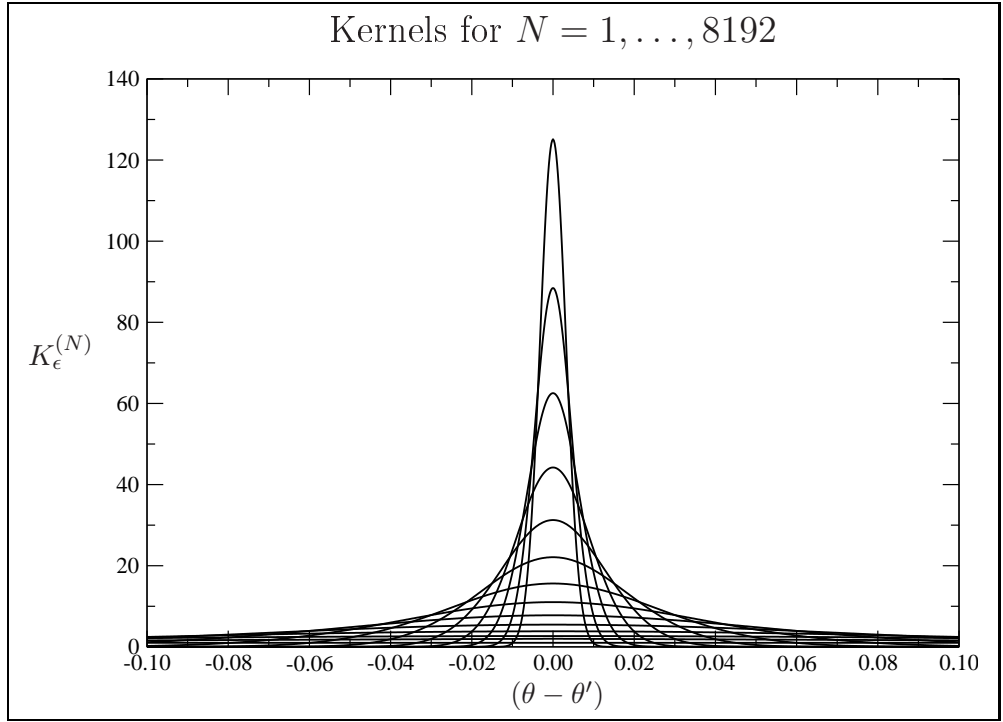,scale=1.0,angle=0}
  }
  \caption{The kernels of several filters with constant range $\epsilon$,
    obtained via the use of their Fourier series, for $\epsilon=0.5$, for
    values of $N$ increasing exponentially, in the set $\{1,2,4,8,16,32,
    64,128,256,512,1024,2048,4096,8192\}$, plotted as functions of
    $\left(\theta-\theta'\right)$ over a small sub-interval within the
    periodic interval $[-\pi,\pi]$.}
  \label{Fig04}
\end{figure}

\noindent
Since this Fourier series converges very fast, and ever faster as $N$
increases, it is very easy to use it to plot the corresponding functions.
Doing this one gets the sequence of functions shown in
Figure~\ref{Fig03}. As expected, the range increases without bound and the
kernel gets distributed more and more over the whole periodic interval,
approaching a constant function with unit integral. This means that only
the constant term of the Fourier expansion survives the limit, and that
all the other Fourier coefficients converge to zero. There is nothing too
surprising about this, since it is consistent with the fact that

\begin{displaymath}
  \lim_{N\to\infty}
  \left[
    \frac{\sin(k\epsilon)}{(k\epsilon)}
  \right]^{N}
  =
  0,
\end{displaymath}

\noindent
so long as $\sin(k\epsilon)<(k\epsilon)$, which is true since $k>0$ and
$\epsilon>0$. Once again a more rigorous analysis of this limit would
require more careful consideration of the convergence of the series, since
in principle one must be careful with the interchange of the $N\to\infty$
limit and the $k\to\infty$ limit of the series. However, here too it turns
out that this rough discussion suffices for our purposes. It is quite
clear that, if we could repeat the experiment on the whole real line
instead of the periodic interval, the kernel would approach a normalized
Gaussian function that in turn would approach zero everywhere, becoming
ever wider and lower as $N\to\infty$.

\begin{figure}[ht]
  \centering
  \fbox{
    \epsfig{file=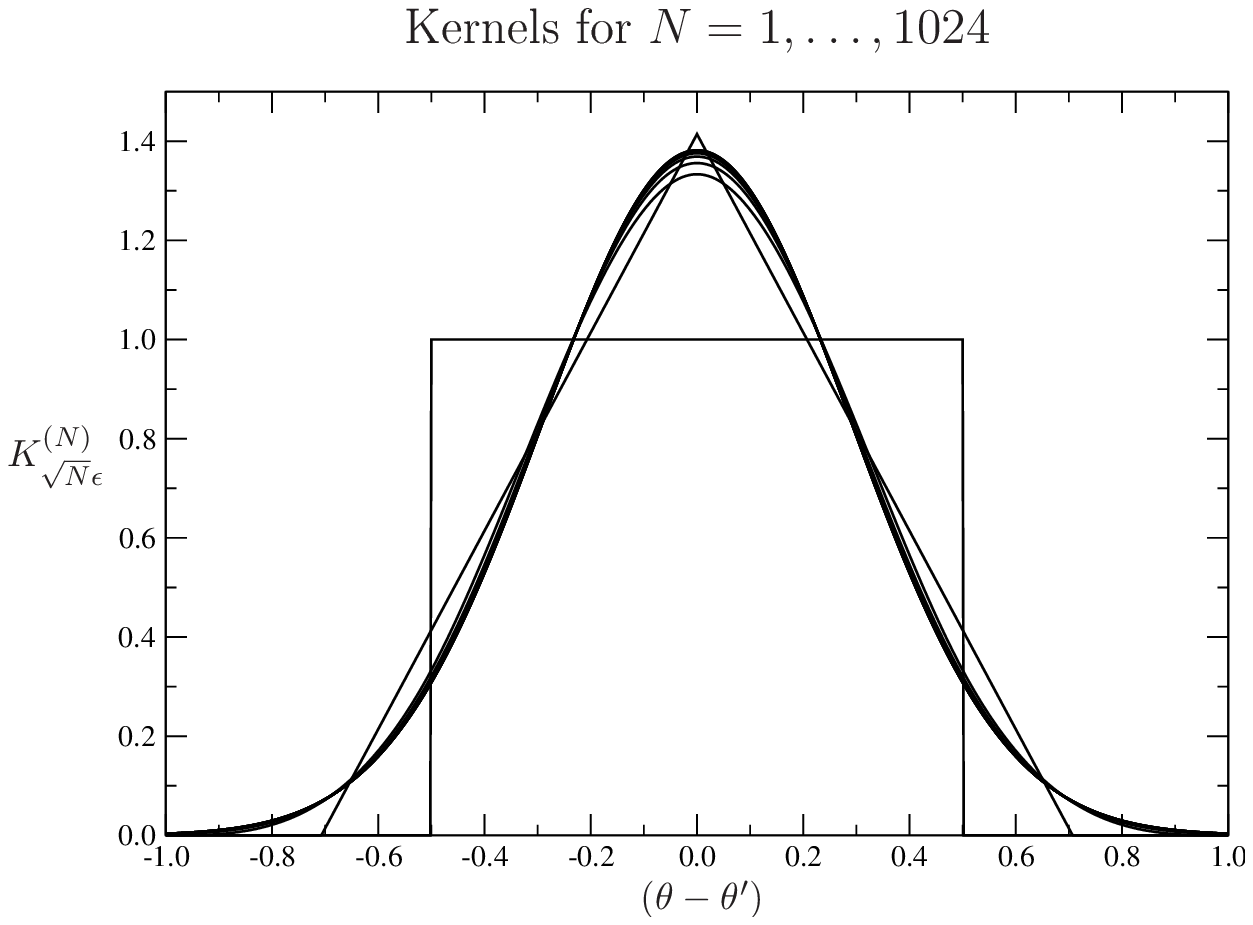,scale=1.0,angle=0}
  }
  \caption{The kernels of several filters with Gaussian range
    $\sqrt{N}\epsilon$, obtained via the use of their Fourier series, for
    $\epsilon=0.5$, for values of $N$ increasing exponentially, in the set
    $\{1,2,4,8,16,32,64,128,256,512,1024\}$, plotted as functions of
    $\left(\theta-\theta'\right)$ over a sub-interval within the periodic
    interval $[-\pi,\pi]$.}
  \label{Fig05}
\end{figure}

If we now consider once again the $N\to\infty$ limit of the kernel of
order $N$ with constant range $\epsilon$, we might try to define an
infinite-order kernel of finite range $\epsilon$ as

\begin{displaymath}
  K_{\epsilon}^{(\infty)}\!\left(\theta-\theta'\right)
  =
  \frac{1}{2\pi}
  +
  \frac{1}{\pi}
  \lim_{N\to\infty}
  \sum_{k=1}^{\infty}
  \left[
    \frac{\sin(k\epsilon/N)}{(k\epsilon/N)}
  \right]^{N}
  \cos\!\left[k\left(\theta-\theta'\right)\right].
\end{displaymath}

\noindent
Our previous analysis indicated that this has the delta ``function'' as
its limit. Plotting this kernel one gets the sequence of functions shown
in Figure~\ref{Fig04}. As one can see, the kernel in fact diverges to
positive infinity at zero. It also seems to go to zero everywhere
else. Since it still has constant integral, and since it can be easily
verified that its maximum at zero diverges to infinity as $\sqrt{N}$, we
must conclude that it in fact approaches a Dirac delta ``function''. More
precisely, the sequence of kernels approaches a normalized Gaussian
function that in turns approaches the delta ``function'', becoming ever
taller and narrower as $N\to\infty$, with constant area under the
graph. What we have here is a singular limit, in a way going full circle,
from the delta ``function'' at $N=0$ and back to it at
$N\to\infty$. Therefore the $N\to\infty$ limit of this order-$N$ kernel is
not a $C^{\infty}$ function, but a singular object instead, which is
certainly an unexpected and surprising result.

This suggests that one may take an intermediate limit, which perhaps will
converge to a non-singular localized function, by superposing $N$ filters
with range $\sqrt{N}\epsilon$, thus obtaining a second infinite-order
kernel with infinite range, given by

\begin{displaymath}
  K_{\infty}^{(\infty)}\!\left(\theta-\theta'\right)
  =
  \frac{1}{2\pi}
  +
  \frac{1}{\pi}
  \lim_{N\to\infty}
  \sum_{k=1}^{\infty}
  \left[
    \frac
    {\sin\!\left(k\epsilon/\sqrt{N}\right)}
    {\left(k\epsilon/\sqrt{N}\right)}
  \right]^{N}
  \cos\!\left[k\left(\theta-\theta'\right)\right].
\end{displaymath}

\noindent
In fact, this exercise results in the sequence of functions seen in
Figure~\ref{Fig05}, which approaches very fast a function very similar to
a normalized Gaussian function with finite and non-zero width. Again it is
quite clear that, if we could repeat this construction on the whole real
line, then the limit would be a normalized Gaussian function, but in our
case here it is a bit deformed by its containment within the periodic
interval. When its width is small when compared to the length $2\pi$ of
the periodic interval the Gaussian approaches zero very fast when we go
significantly away from its point of maximum, so that we may consider this
case to be paracompact, and even use this filter successfully in the
practice of physics applications. However, the exact mathematical fact is
that the range of this filter does tend to infinity when $N\to\infty$, and
therefore to the whole extent of the periodic interval, if we execute all
this operation within it.

It is possible to interpret qualitatively what happens in these three
cases in terms of the expression for the corresponding superpositions of
inner analytic functions. In the case of the straight multiple
superposition of the first-order kernel we have for the representation of
the order-$N$ filter on the complex plane

\begin{displaymath}
  w_{N\epsilon}^{(N)}(z)
  =
  \left(\frac{-\ii}{2\epsilon}\right)^{N}
  \sum_{n=0}^{N}
  \frac{(-1)^{n}N!}{n!(N-n)!}\,
  w^{-N\ldot}\!\left(z\e{\ii[N-2n]\epsilon}\right),
\end{displaymath}

\noindent
where we may assume that the coefficients do not diverge with $N$, since
the function tends to a constant everywhere on the unit circle in the
$N\to\infty$ limit. In the case of the superposition of the first-order
kernel with decreasing range $\epsilon/N$, resulting on a fixed range
$\epsilon$ for the order-$N$ kernel, the representation of the order-$N$
filter in the complex plane is

\begin{displaymath}
  w_{\epsilon}^{(N)}(z)
  =
  \left(\frac{-\ii N}{2\epsilon}\right)^{N}
  \sum_{n=0}^{N}
  \frac{(-1)^{n}N!}{n!(N-n)!}\,
  w^{-N\ldot}\!\left(z\e{\ii[1-2n/N]\epsilon}\right),
\end{displaymath}

\noindent
so that the extra factor of $N^{N}$ is clearly related to the divergence
at a point on the unit circle, when we make $N\to\infty$. In the case of
the superposition of the first-order kernel with decreasing range
$\epsilon/\sqrt{N}$, resulting on a range $\sqrt{N}\epsilon$ for the
order-$N$ kernel, the representation of the order-$N$ filter in the
complex plane is

\begin{displaymath}
  w_{\sqrt{N}\epsilon}^{(N)}(z)
  =
  \left(\frac{-\ii\sqrt{N}}{2\epsilon}\right)^{N}
  \sum_{n=0}^{N}
  \frac{(-1)^{n}N!}{n!(N-n)!}\,
  w^{-N\ldot}\!
  \left(z\e{\ii\left[\sqrt{N}-2n/\sqrt{N}\right]\epsilon}\right).
\end{displaymath}

\noindent
Note that in this case we gained a factor of $N^{N/2}$ rather than
$N^{N}$, which the consequence that over the unit circle the kernel
neither approaches a constant everywhere nor diverges to infinity
somewhere. In this case the coefficients seem to have well-defined finite
limits.

We must therefore conclude that, with this type of multiple superposition
of the first-order filter, and the corresponding superposition of the
singularities of the inner analytic functions in the complex plane, we are
unable to define an infinite-order kernel that is both a finite and smooth
real function, and that at the same time is localized within a compact
support, thus generating an infinite-order filter with a finite range
$\epsilon$. In order to understand why, it is useful to look at the
singularities, on the unit circle of the complex plane, of the sequence of
inner analytic functions generated by the repeated application of the
first-order filter, starting with the inner analytic function
corresponding to the zero-order kernel, which is a Dirac delta
``function'' and thus has a single first-order pole at some point on the
unit circle, as shown in~\cite{FTotCPI}.

If we look at the diagram in Figure~\ref{Fig02}, we see that as the
multiple application of the first-order filter goes on, more and more
softened singularities are superposed at the points near the center of the
diagram. Each singularity is progressively softer, but they are superposed
in increasing numbers, thus generating a coefficient in the corresponding
term in the superposition shown in Equation~(\ref{filterordNrangee}). For
each finite $N$ these coefficients may be large, but they are finite, and
therefore they do not disturb the softness of the corresponding
singularities. However, if one of the coefficients diverges in the
$N\to\infty$ limit, then the corresponding singularity is no longer soft
in the limit. Let us recall that the definition of a soft singularity, as
given in~\cite{FTotCPII}, is that the limit of the inner analytic function
to that point be finite. Because of the diverging coefficients, in this
type of superposition this may fail to be so in the $N\to\infty$ limit,
even if the singularities are soft for each finite value of $N$.

\begin{figure}[ht]
  \centering
  \fbox{
    \epsfig{file=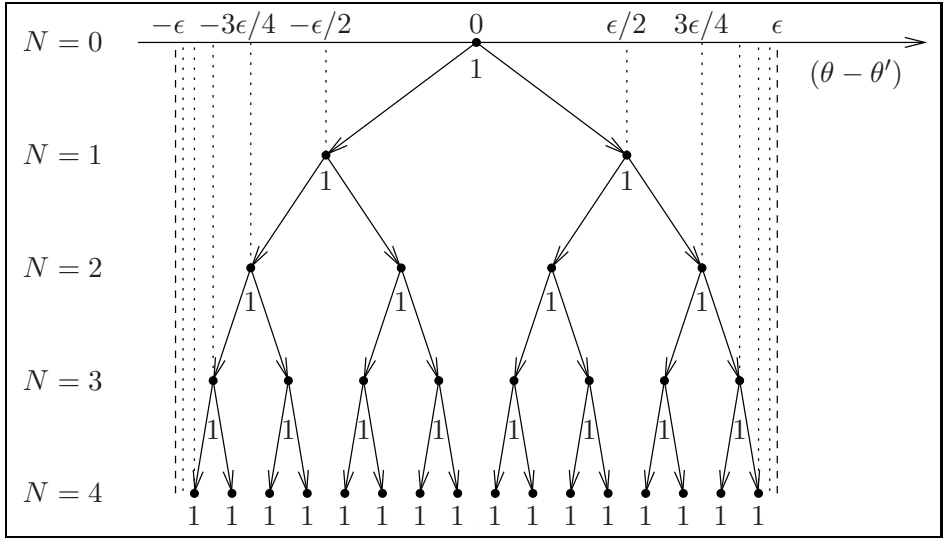,scale=1.0,angle=0}
  }
  \caption{The scaled iteration of the first-order filter to produce an
    order-$N$ filter, showing the range $\epsilon_{N}$ tending to the
    limit $\epsilon$. The original singularity is at $\theta'$. The
    numbers near the vertices of the triangles show that there is just one
    softened singularity at each such point.}
  \label{Fig06}
\end{figure}

We are therefore led to the idea of changing the method of iteration of
the first-order filters in such a way that the softened singularities
never get superposed. It is not too difficult to see that one may
accomplish this by superposing filters with progressively smaller ranges,
as illustrated by the diagram in Figure~\ref{Fig06}. Given a value of
$\epsilon$, this diagram corresponds to a process in which we start by
applying the first-order filter with range $\epsilon/2$, followed by the
first-order filter of range $\epsilon/4$, then by the filter of range
$\epsilon/8$, and so on, where the range of the $N^{\rm th}$ filter
applied is given by $\epsilon/2^{N}$. Note that at the $N^{\rm th}$
iteration there are $2^{N}$ singularities homogeneously distributed within
the interval $(-\epsilon,\epsilon)$. Since the ranges are scaled down
exponentially, we call this a {\em scaled filter}, and the corresponding
kernel a {\em scaled kernel}.

\begin{figure}[ht]
  \centering
  \fbox{
    \epsfig{file=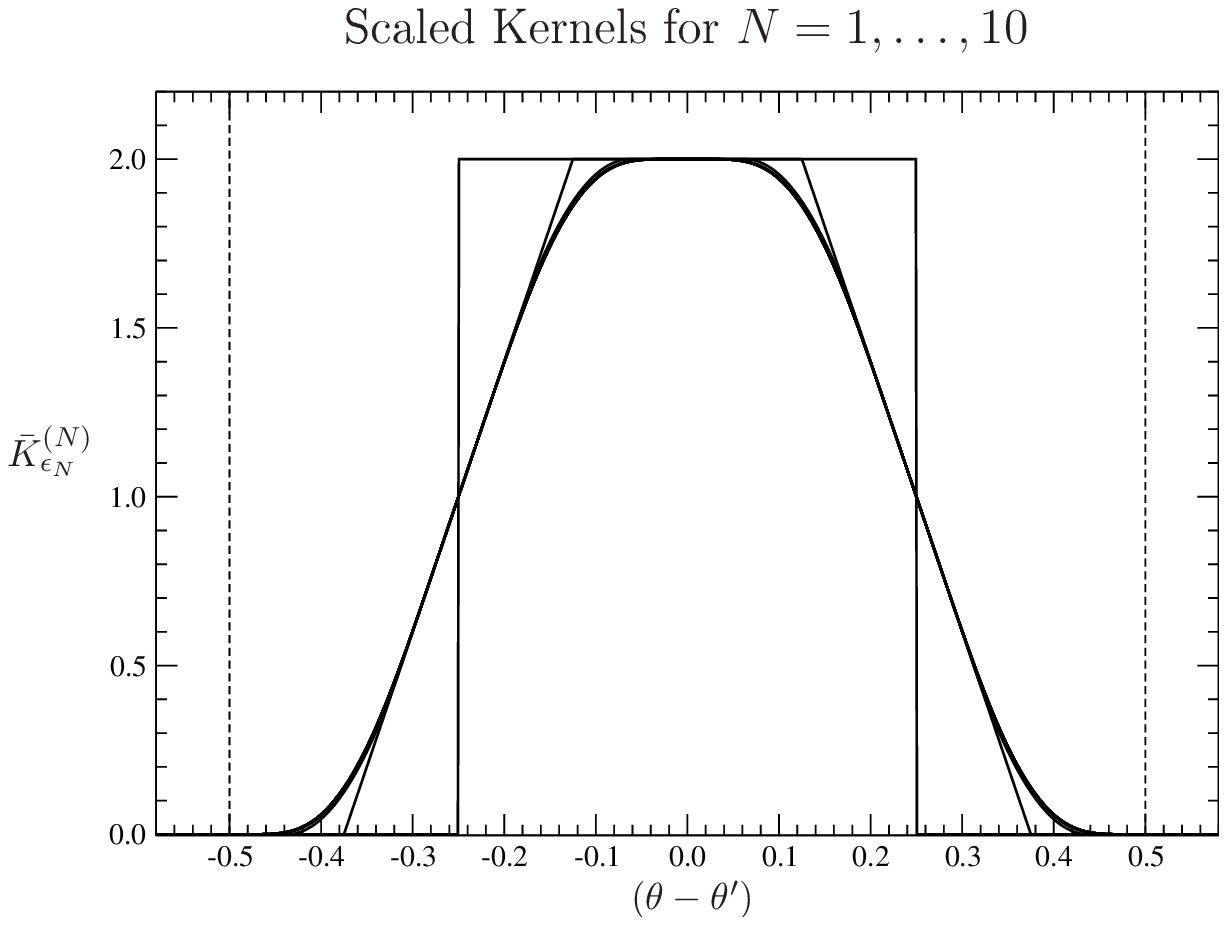,scale=1.0,angle=0}
  }
  \caption{The kernels of several scaled filters with range
    $\epsilon_{N}\to\epsilon$, obtained via the use of their Fourier
    series, for $\epsilon=0.5$, for values of $N$ increasing linearly, in
    the set $\{1,2,3,4,5,6,7,8,9,10\}$, plotted as functions of
    $\left(\theta-\theta'\right)$ over the support interval
    $[-\epsilon,\epsilon]$. The dashed lines mark the ends of the support
    interval.}
  \label{Fig07}
\end{figure}

At the $N^{\rm th}$ iteration there are $2^{N}$ singularities, each one
softened by $N$ degrees, spaced from one another by $\epsilon/2^{N-1}$,
and spaced from the ends of the $[-\epsilon,\epsilon]$ interval by
$\epsilon/2^{N}$. They are, therefore, regularly distributed within the
interval $(-\epsilon,\epsilon)$, centered at $2^{N}$ consecutive
sub-intervals of length $\epsilon/2^{N-1}$. In the $N\to\infty$ limit the
singularities will tend to become homogeneously distributed within the
interval $(-\epsilon,\epsilon)$, and indeed will tend to a countable
infinity of infinitely soft singularities distributed densely within that
interval. Since at the $N^{\rm th}$ step there are $2^{N}$ such
singularities, we see that their number grows exponentially fast. However,
the singularities never get superposed. The corresponding superposition in
terms of inner analytic functions in the complex plane is given by

\noindent
\begin{eqnarray*}
  \bar{w}_{\epsilon_{N}}^{(N)}(z)
  & = &
  \left(\frac{-2^{1}\ii}{\epsilon}\right)
  \left(\frac{-2^{2}\ii}{\epsilon}\right)
  \times\ldots\times
  \left(\frac{-2^{N-1}\ii}{\epsilon}\right)
  \left(\frac{-2^{N}\ii}{\epsilon}\right)
  \times
  \\
  &   &
  \hspace{7em}
  \times
  \sum_{n=1}^{2^{N}}
  (-1)^{n-1}
  w^{-N\ldot}\!\left(z\e{\ii\left[1-(2n-1)/2^{N}\right]\epsilon}\right)
  \\
  & = &
  \left(\frac{-\ii}{\epsilon}\right)^{N}
  2^{N(N+1)/2}
  \sum_{n=1}^{2^{N}}
  (-1)^{n-1}
  w^{-N\ldot}\!\left(z\e{\ii\left[1-(2n-1)/2^{N}\right]\epsilon}\right).
\end{eqnarray*}

\noindent
Using the property of the first-order filter regarding its action on
Fourier expansions~\cite{LPFprop11,LPFprop12}, it is not difficult to
write the Fourier expansion of this new scaled kernel, at the $N^{\rm th}$
step of the construction process, which has a range $\epsilon_{N}$ such
that $\epsilon/2\leq \epsilon_{N}<\epsilon$,

\begin{displaymath}
  \bar{K}_{\epsilon_{N}}^{(N)}\!\left(\theta-\theta'\right)
  =
  \frac{1}{2\pi}
  +
  \frac{1}{\pi}
  \sum_{k=1}^{\infty}
  \left[
    \frac
    {\sin\!\left(k\epsilon/2^{1}\right)}
    {\left(k\epsilon/2^{1}\right)}
  \right]
  \times
  \ldots
  \times
  \left[
    \frac
    {\sin\!\left(k\epsilon/2^{N}\right)}
    {\left(k\epsilon/2^{N}\right)}
  \right]
  \cos\!\left[k\left(\theta-\theta'\right)\right],
\end{displaymath}

\noindent
where the coefficients include a product of $N$ different sinc factors.
This product can be written as

\noindent
\begin{eqnarray*}
  \left[
    \frac
    {\sin\!\left(k\epsilon/2^{1}\right)}
    {\left(k\epsilon/2^{1}\right)}
  \right]
  \times
  \ldots
  \times
  \left[
    \frac
    {\sin\!\left(k\epsilon/2^{N}\right)}
    {\left(k\epsilon/2^{N}\right)}
  \right]
  & = &
  \frac{2^{1+\ldots+N}}{(k\epsilon)^{N}}\,
  \sin\!\left(\frac{k\epsilon}{2^{1}}\right)
  \times
  \ldots
  \times
  \sin\!\left(\frac{k\epsilon}{2^{N}}\right)
  \\
  & = &
  \frac{2^{N(N+1)/2}}{(k\epsilon)^{N}}
  \prod_{n=1}^{N}
  \sin\!\left(\frac{k\epsilon}{2^{n}}\right).
\end{eqnarray*}

\noindent
It follows that we have for this order-$N$ scaled kernel

\begin{equation}\label{expanfilt}
  \bar{K}_{\epsilon_{N}}^{(N)}\!\left(\theta-\theta'\right)
  =
  \frac{1}{2\pi}
  +
  \frac{1}{\pi}
  \sum_{k=1}^{\infty}
  \frac{2^{N(N+1)/2}}{(k\epsilon)^{N}}
  \left[
    \,
    \prod_{n=1}^{N}
    \sin\!\left(\frac{k\epsilon}{2^{n}}\right)
  \right]
  \cos\!\left[k\left(\theta-\theta'\right)\right].
\end{equation}

\noindent
As expected, this kernel indeed has a well-defined limit when
$N\to\infty$, with support in the finite interval
$\left[\theta'-\epsilon,\theta'+\epsilon\right]$, as one can see in the
graph of Figure~\ref{Fig07}. The sequence of kernels converges
exponentially fast to a definite function, with a rather unusual shape.

\begin{figure}[ht]
  \centering
  \fbox{
    \epsfig{file=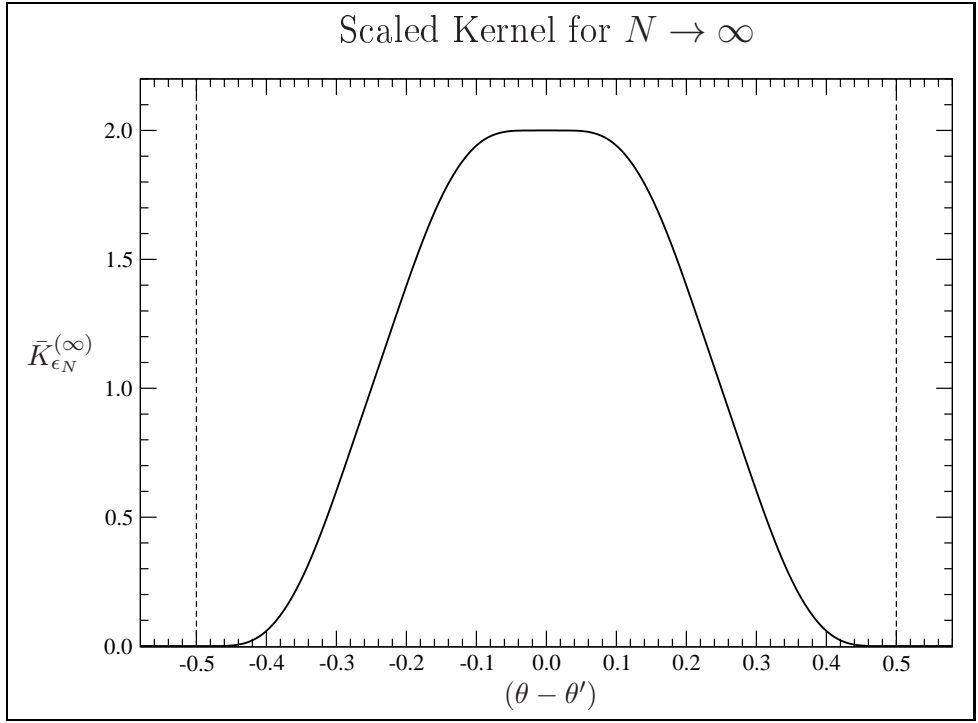,scale=1.0,angle=0}
  }
  \caption{The best approximation of the kernel of the scaled filter with
    constant range $\epsilon$, in the $N\to\infty$ limit, obtained via the
    use of its Fourier series, for $\epsilon=0.5$, for a large value of
    $N$ ($100$), plotted as a function of $\left(\theta-\theta'\right)$
    over the support interval $[-\epsilon,\epsilon]$. The dashed lines
    mark the ends of the support interval.}
  \label{Fig08}
\end{figure}

It is possible to demonstrate explicitly the convergence of the sequence
of scaled kernels
$\bar{K}_{\epsilon_{N}}^{(N)}\!\left(\theta-\theta'\right)$ to a
well-defined regular function
$\bar{K}_{\epsilon}^{(\infty)}\!\left(\theta-\theta'\right)$ in the
$N\to\infty$ limit. The proof is rather lengthy and is presented in full
in Appendix~\ref{APPInfOrd}. It depends on the following facts about this
limit, that we may establish here in order to give a general idea of the
structure of the proof. First of all, due to one of the properties of the
first-order filter~\cite{LPFprop09} all kernels in the sequence have unit
integral. Second, the range $\epsilon_{N}$ of the order-$N$ kernel is
given by the combined ranges of all the kernels used to build it, and is
therefore given by

\begin{displaymath}
  \epsilon_{N}
  =
  \frac{\epsilon}{2}
  +
  \frac{\epsilon}{4}
  +
  \ldots
  +
  \frac{\epsilon}{2^{N-1}}
  +
  \frac{\epsilon}{2^{N}},
\end{displaymath}

\noindent
which is a geometric progression with ratio $1/2$. We have therefore

\noindent
\begin{eqnarray*}
  \epsilon_{N}
  & = &
  \epsilon\,
  \frac
  {
    \FFrac{1}{2}-\FFrac{1}{2^{N+1}}
  }
  {
    1-\FFrac{1}{2}
  }
  \\
  & = &
  \epsilon
  \left(
    1-\frac{1}{2^{N}}
  \right).
\end{eqnarray*}

\noindent
It follows therefore that in the $N\to\infty$ limit $\epsilon_{N}$ tends
to $\epsilon$,

\begin{displaymath}
  \lim_{N\to\infty}\epsilon_{N}
  =
  \epsilon.
\end{displaymath}

\noindent
Therefore all the kernels in the sequence remain identically zero
everywhere outside the interval $(-\epsilon,\epsilon)$, so that we may
conclude that the limiting function has support within that interval. Next
we observe that, since the filtered function is defined as an average of
the original function, it can never assume values which are larger than
the maximum of the function it is applied on, or smaller than its minimum.
Therefore, since the first kernel we start with, with range $\epsilon/2$,
is limited within the interval $[0,1/\epsilon]$, so are all the subsequent
kernels of the construction sequence.

As a consequence all these considerations, the non-zero parts of the
graphs of all the kernels in the construction sequence are contained
within the rectangle defined by the support interval
$[-\epsilon,\epsilon]$ and the range of values $[0,1/\epsilon]$, which has
area $2$. The area of the graph of every kernel in the construction
sequence is $1$, so that it occupies one half of the area of the
rectangle. Since every scaled kernel in the construction sequence is a
continuous and differentiable function, it becomes clear that the
$N\to\infty$ limit of the sequence of kernels must also be a regular
function within this rectangle. It follows from the discussion in
Appendix~\ref{APPInfOrd} that the limit of the sequence exists and is a
regular, continuous and differentiable function, resulting in the
infinite-order scaled kernel with range $\epsilon$

\begin{displaymath}
  \bar{K}_{\epsilon}^{(\infty)}\!\left(\theta-\theta'\right)
  =
  \lim_{N\to\infty}
  \bar{K}_{\epsilon_{N}}^{(N)}\!\left(\theta-\theta'\right).
\end{displaymath}

\noindent
This infinite-order scaled kernel has the deceptively simple look shown in
Figure~\ref{Fig08}. Despite appearances it contains no completely straight
segments within its support interval. Once it is shown that it is a
$C^{\infty}$ function, it follows that its derivatives of all orders are
zero at the two extremes of the support interval, where they must match
the correspondingly zero derivatives of the two external segments, on
either side of the support interval, where the kernel is identically zero.

A technical note about the graphs representing the $N\to\infty$ limit of
various quantities is in order at this point. They were obtained
numerically from the corresponding Fourier series, using the scaled filter
of order $N=100$. This means that the last first-order filter used in the
multiple superposition has a range of $\epsilon/2^{100}$. This is less
than $\epsilon/10^{30}$ and is therefore many orders of magnitude below
any graphical resolution one might hope for, in any medium. Certainly the
errors related to the summation of the Fourier series, which were set at
$10^{-12}$, are the dominant ones, but still extremely small. We may
conclude that these graphs are faithful representations of the
corresponding quantities for all conceivable graphical purposes. The
programs used in creating all the graphs shown in this paper are freely
available online~\cite{TarFile}.

It is quite simple to see that this kernel is a $C^{\infty}$ function. Its
Fourier series, given as the $N\to\infty$ limit of the Fourier expansion
in Equation~(\ref{expanfilt}), is certainly absolutely and uniformly
convergent, and any finite-order term-wise derivative of it results in
another series with the same properties. Since the convergence of the
resulting series is the additional condition, besides uniform convergence,
that suffices to guarantee that one can differentiate the series term-wise
in order to obtain the derivative of the function, we may conclude that
derivatives of all finite orders exist and are given by continuous and
differentiable functions. This proof of infinite differentiability of the
$\bar{K}_{\epsilon}^{(\infty)}\!\left(\theta-\theta'\right)$ kernel now
ensures that the kernel and all its multiple derivatives are in fact zero
at the points $\left(\theta-\theta'\right)=\pm\epsilon$.

An independent discussion of the existence of the derivatives of all
finite orders can be found in Section~\ref{APPInfOrdDiff} of
Appendix~\ref{APPInfOrd}. There one can see also that the derivatives of
all orders are zero at the central point of maximum, as well as at the two
extremes of the support interval. We also show in
Section~\ref{APPNonAnalytic} of Appendix~\ref{APPInfOrd} that the
infinite-order scaled kernel is {\em not} analytic as a function on the
periodic interval, in the real sense of the term. We do this by showing
that there is a countable infinity of points, distributed densely in the
support interval, where only a finite number of derivatives is different
from zero. We also discuss briefly there the question of whether or not
the infinite-order scaled kernel can be extended analytically to the
complex plane. We discuss this in terms of the fact that it is the limit
of an inner analytic function when one takes the limit to the border of
the unit disk.

\begin{figure}[ht]
  \centering
  \fbox{
    \epsfig{file=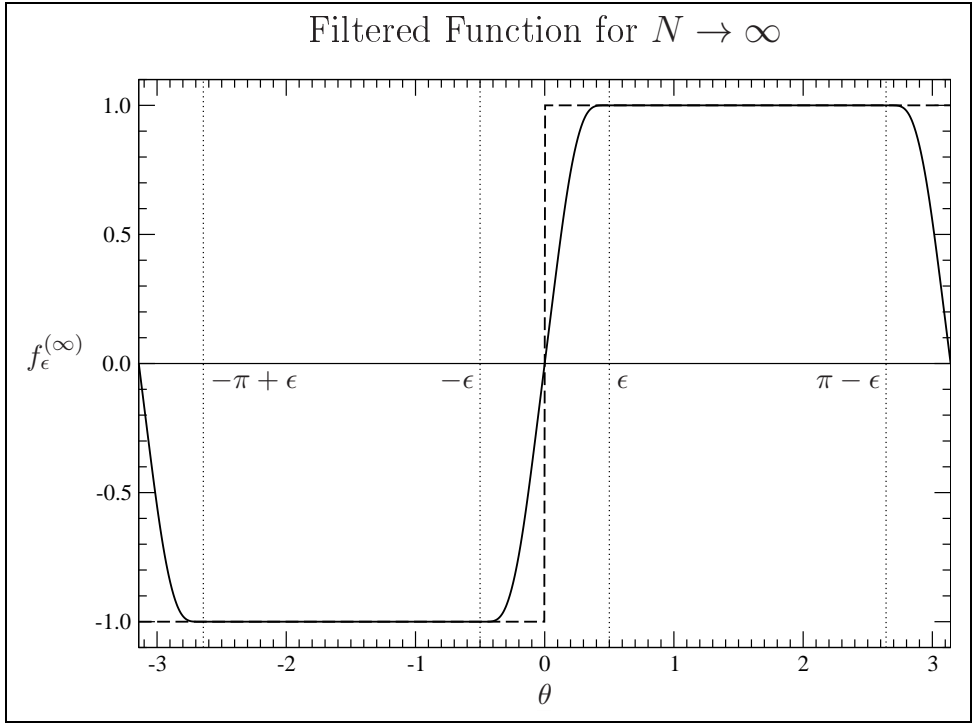,scale=1.0,angle=0}
  }
  \caption{The filtered square wave and the parameters related to the
    action of the infinite-order scaled filter, obtained via the use of
    its Fourier series, for $\epsilon=0.5$, for a large value of $N$
    ($100$), plotted as a function of $\theta$ over the periodic interval
    $[-\pi,\pi]$. The original function is shown with the dashed line and
    the filtered function with the solid line. The dotted lines mark the
    intervals where the function was changed by the filter.}
  \label{Fig09}
\end{figure}

Given that the infinite-order scaled kernel is a $C^{\infty}$ function,
one may then define an infinite-order filter based on this infinite-order
scaled kernel. However, in this case it is not so simple to determine the
form of the coefficients directly in the $N\to\infty$ limit. In fact, the
question of whether or not it is possible to write the coefficients
directly in the $N\to\infty$ limit, in some simple form, is an open
one. Note that, as was mentioned before, one must be careful with the
interchange of the $N\to\infty$ limit and the $k\to\infty$ limit of the
series. In any case, it follows that given {\em any} merely integrable
function $f(\theta)$, the filtered function

\begin{displaymath}
  f_{\epsilon}^{(\infty)}(\theta)
  =
  \int_{-\infty}^{\infty}d\theta'\,
  \bar{K}_{\epsilon}^{(\infty)}\!\left(\theta-\theta'\right)
  f\!\left(\theta'\right),
\end{displaymath}

\noindent
is necessarily a $C^{\infty}$ function, in the real sense of the term. It
is easy to see this, since the differentiations with respect to $\theta$
at the right-hand side will act only on the infinite-order scaled kernel,
which is a $C^{\infty}$ function. Therefore all the finite-order
derivatives of $f_{\epsilon}^{(\infty)}(\theta)$ exist, since they may be
written as

\begin{displaymath}
  \frac{d^{n}}{d\theta^{n}}
  f_{\epsilon}^{(\infty)}(\theta)
  =
  \int_{-\infty}^{\infty}d\theta'
  \left[
    \frac{d^{n}}{d\theta^{n}}
    \bar{K}_{\epsilon}^{(\infty)}\!\left(\theta-\theta'\right)
  \right]
  f\!\left(\theta'\right),
\end{displaymath}

\noindent
where the $n^{\rm th}$ derivative of the infinite-order scaled kernel is a
limited, continuous and differentiable function with compact support,
being therefore an integrable function. In addition to this, the changes
made in $f(\theta)$ in order to produce $f_{\epsilon}^{(\infty)}(\theta)$
have a finite range $\epsilon$ that can be made as small as one wishes.

\begin{figure}[ht]
  \centering
  \fbox{
    \epsfig{file=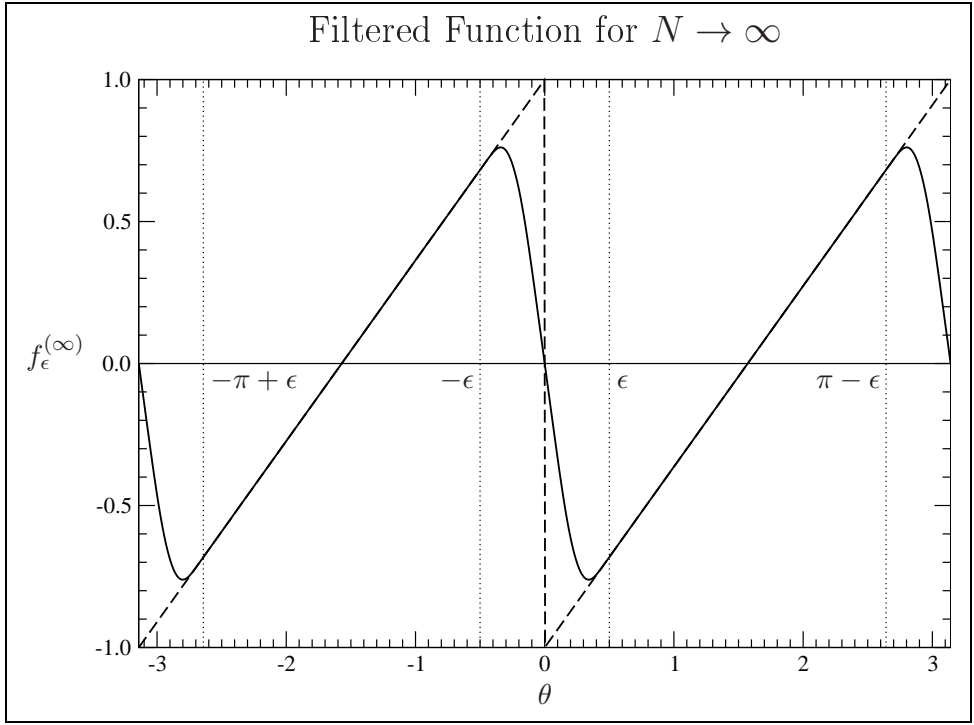,scale=1.0,angle=0}
  }
  \caption{The filtered sawtooth wave wave and the parameters related to
    the action of the infinite-order scaled filter, obtained via the use
    of its Fourier series, for $\epsilon=0.5$, for a large value of $N$
    ($100$), plotted as a function of $\theta$ over the periodic interval
    $[-\pi,\pi]$. The original function is shown with the dashed line and
    the filtered function with the solid line. The dotted lines mark the
    intervals where the function was changed by the filter.}
  \label{Fig10}
\end{figure}

Let us now consider the action of such a filter on inner analytic
functions. There is no difficulty in determining what happens to the
singularities of the inner analytic functions on the unit circle, as a
consequence of the application of this infinite-order scaled filter. It is
quite clear that a single singularity of the inner analytic function at
$\theta$ would be smeared into a denumerable infinity of softened
singularities within the interval $[\theta-\epsilon,\theta+\epsilon]$ on
the unit circle. This set of singularities would occupy the interval
densely, and they would be infinitely soft, since they are the result of
an infinite sequence of logarithmic integrations. Since each logarithmic
integration renders the real function on the unit circle differentiable to
one more order, in the limit one gets over that circle an infinitely
differentiable real function of $\theta$. So the resulting complex
function of $z$ must be an inner analytic function which has only
infinitely soft singularities on the unit circle and that is a
$C^{\infty}$ function of $\theta$ when restricted to that circle. The same
is true for the kernel itself, if we start with a single first-order pole
at $\theta$, which is the case for the order-zero kernel. Note that in
either case this real function is $C^{\infty}$ right on top of a
densely-distributed set of singularities of the corresponding inner
analytic function, which is somewhat unexpected, even if they are
infinitely soft singularities.

As a simple example, let us consider the unit-amplitude square wave, which
is a discontinuous periodic function, as one can see in
Figure~\ref{Fig09}, that shows the filtered function superposed with the
original one. The filtered function was obtained from its Fourier series,
which due to the properties of the first-order filter is easily obtained,
being given by

\begin{displaymath}
  f_{\epsilon_{N}}(\theta)
  =
  \frac{4}{\pi}
  \sum_{j=0}^{\infty}
  \frac{2^{N(N+1)/2}}{k^{N+1}\epsilon^{N}}
  \left[
    \,
    \prod_{n=1}^{N}
    \sin\!\left(\frac{k\epsilon}{2^{n}}\right)
  \right]
  \sin\!\left(k\theta\right),
\end{displaymath}

\noindent
where $k=2j+1$, for a large value of $N$. The graph of the original
function has two straight horizontal segments and two points of
discontinuity at $\theta=0$ and at $\theta=\pm\pi$. It follows that the
corresponding inner analytic function has two borderline hard
singularities at these two points. Let us consider all the instances of
the first-order linear low-pass filter used for the construction of the
infinite-order scaled kernel, for all construction steps $N$ and any value
of $\epsilon<\pi$. Since the linear low-pass filters are all the identity
on the segments that are linear functions, up to a distance of $\epsilon$
to one of the singularities, the function would never be changed at all
outside the two intervals $[-\epsilon,\epsilon]$ and
$[\pi-\epsilon,-\pi+\epsilon]$, when one applies to it any of the
order-$N$ scaled filters. After the end of the process of application of
the infinite-order scaled filter these two intervals would contain
segments of $C^{\infty}$ functions of $\theta$, and in fact the whole
resulting function would be a $C^{\infty}$ function of $\theta$, over the
whole unit circle.

\begin{figure}[ht]
  \centering
  \fbox{
    \epsfig{file=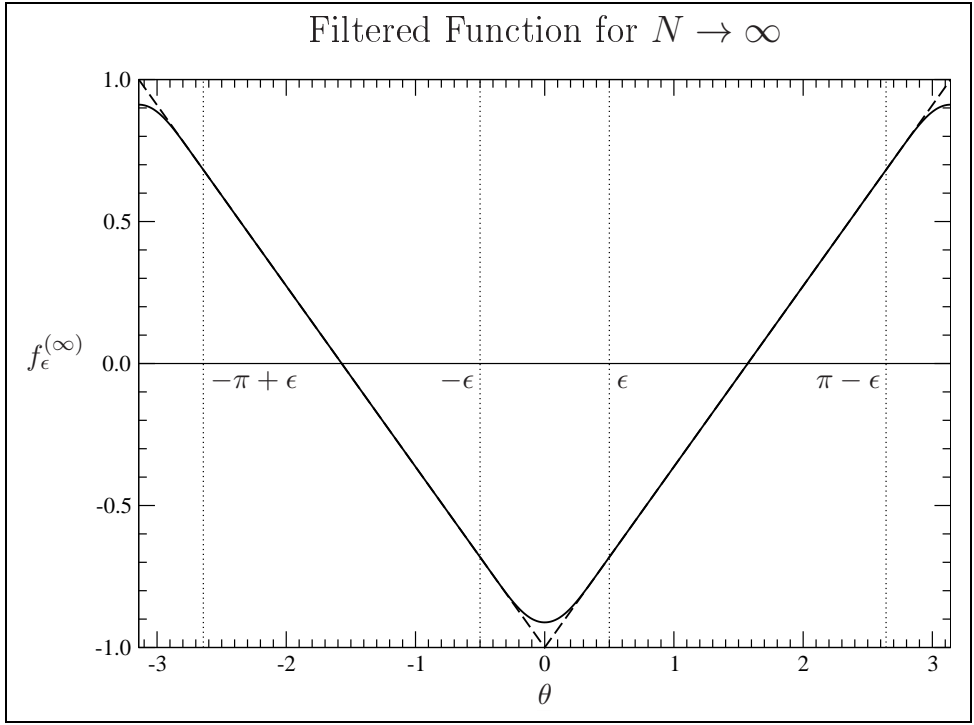,scale=1.0,angle=0}
  }
  \caption{The filtered triangular wave and the parameters related to the
    action of the infinite-order scaled filter, obtained via the use of
    its Fourier series, for $\epsilon=0.5$, for a large value of $N$
    ($100$), plotted as a function of $\theta$ over the periodic interval
    $[-\pi,\pi]$. The original function is shown with the dashed line and
    the filtered function with the solid line. The dotted lines mark the
    intervals where the function was changed by the filter.}
  \label{Fig11}
\end{figure}

Another similar example can be seen in Figure~\ref{Fig10}, which shows the
case of the unit-amplitude sawtooth wave, which also has the same two
points of discontinuity and therefore corresponds to an inner analytic
function with two similar borderline hard singularities at these points.
The filtered function was obtained from its Fourier series,

\begin{displaymath}
  f_{\epsilon_{N}}(\theta)
  =
  -\,
  \frac{4}{\pi}
  \sum_{j=1}^{\infty}
  \frac{2^{N(N+1)/2}}{k^{N+1}\epsilon^{N}}
  \left[
    \,
    \prod_{n=1}^{N}
    \sin\!\left(\frac{k\epsilon}{2^{n}}\right)
  \right]
  \sin\!\left(k\theta\right),
\end{displaymath}

\noindent
where $k=2j$, for a large value of $N$. In Figure~\ref{Fig11} one can see
the case of the triangular wave, which is a continuous function with two
points of non-differentiability at $\theta=0$ and at $\theta=\pm\pi$, and
therefore corresponds to an inner analytic function with two borderline
soft singularities at these points. The filtered function was obtained
from its Fourier series,

\begin{displaymath}
  f_{\epsilon_{N}}(\theta)
  =
  -\,
  \frac{8}{\pi^{2}}
  \sum_{j=0}^{\infty}
  \frac{2^{N(N+1)/2}}{k^{N+2}\epsilon^{N}}
  \left[
    \,
    \prod_{n=1}^{N}
    \sin\!\left(\frac{k\epsilon}{2^{n}}\right)
  \right]
  \cos\!\left(k\theta\right),
\end{displaymath}

\noindent
where $k=2j+1$, for a large value of $N$. In all cases we chose a rather
large value for $\epsilon$ as compared to its maximum value $\pi$, namely
$0.5$, in order to render the action of the scaled infinite-order filter
clearly visible. What we seem to have here is a factory of $C^{\infty}$
functions of $\theta$ on the unit circle. Starting with virtually {\em any
  integrable function}, we may consider the application of the
infinite-order scaled filter in order to produce a $C^{\infty}$ function
on the unit circle, making changes only with a range $\epsilon$ that can
be as small as we wish.

Once we have the infinite-order scaled filter defined within the periodic
interval, it is simple to extend it to the whole real line. Considering
that the infinite-order scaled kernel and all its derivatives are zero at
the two ends of its support interval, we may just take that support
interval and insert it into the real line. If we make the new
infinite-order scaled kernel identically zero outside the support
interval, in all the rest of the real line, we still have a $C^{\infty}$
function. This is so because at the two points of concatenation the two
lateral limits of the kernel are equal, being both zero, as are the two
lateral limits of its first derivative, and the same for all the
higher-order derivatives. Therefore, we may also define an infinite-order
scaled filter acting on the whole real line, that maps any integrable real
function to corresponding $C^{\infty}$ functions.

\section{Conclusions}

Linear low-pass filters of arbitrary orders can be easily and elegantly
defined on the complex plane, within the unit disk, acting on inner
analytic functions. Within the open unit disk the filter simply maps inner
analytic functions onto other inner analytic functions. Through the
correspondence of these functions with FC pairs of DP Fourier series,
these filters reproduce the linear low-pass filters that were defined in a
previous paper, acting on the corresponding DP real functions defined on
the unit circle. Several of the properties of these filters are then
clearly in view, given the known properties of that correspondence.

The effect of the first-order low-pass filter, as seen in the complex
plane, is characterized as a process of singularity splitting, in which
each singularity of an inner analytic function on the unit circle is
exchanged for two softer singularities over that same circle. This has the
effect of improving the convergence characteristics of the DP Fourier
series, and also of rendering the corresponding DP real functions smoother
after the filtering process. Higher-order filters correspond to the
iteration of this process on the unit circle, producing ever larger
collections of ever softer singularities on that circle.

A discussion of the problems encountered when one tries to define an
infinite-order low-pass filter acting on real functions, in the most
immediate way, led to the detailed construction of such an infinite-order
filter, within a compact support. The representation of the filters in the
complex plane was instrumental for the success of this construction. The
infinite-order filter is defined in terms of an infinite-order scaled
kernel, with compact support given by a real parameter $\epsilon$, which
can be as small as one wishes.

The infinite-order scaled kernel is defined as the limit of a sequence of
order-$N$ scaled kernels, and proof of the convergence of the sequence was
presented. It was also shown that the infinite-order scaled kernel is a
$C^{\infty}$ real function, but not an analytic real function. The
infinite-order scaled kernel can be given in a fairly explicit way as a
limit of a Fourier series, which converges extremely fast. Once the
infinite-order filter is defined in the periodic interval, it is a simple
matter to define a corresponding infinite-order filter that acts on the
whole real line.

This infinite-order scaled filter, acting on any merely integrable real
function on the unit circle, has as its result a real function that is
$C^{\infty}$ on the unit circle, while making on the original function
only changers with the finite range $\epsilon$. The same is true for real
function defined on the whole real line. Therefore, one obtains as a
result of this construction a tool that can produce from any integrable
real function corresponding $C^{\infty}$ functions, making changes only
within a range $\epsilon$ that can be as small as desired.

This allows us to use these filters in physics applications, if we use
values of $\epsilon$ sufficiently small in order not to change the
description of the physics within the physically relevant scales of any
given problem. By reducing the value of $\epsilon$ these $C^{\infty}$
functions can be made as close as one wishes to the corresponding original
functions, according to a criterion that has a clear physical meaning, as
explained in a previous paper.

In addition to this, the construction of the filter is equivalent to proof
that there are many real functions that are $C^{\infty}$ but that are not
analytic, and that are typically not extensible analytically to the
complex plane. The filter can be used to produce examples of such
functions in copious quantities. It is quite easy to obtain accurate
values for the filtered functions by numerical means, and thus to
represent the action of the filter in practical applications.

\section{Acknowledgements}

The author would like to thank his friend and colleague Prof. Carlos
Eugênio Imbassay Carneiro, to whom he is deeply indebted for all his
interest and help, as well as his careful reading of the manuscript and
helpful criticism regarding this work.

\appendix

\section{Appendix: Technical Proofs}\label{APPProofs}

\subsection{Direct Derivation of the Coefficients {\boldmath
    $a_{\epsilon,k}$}}\label{APPFourCoef}

Let us determine the effect of the first-order linear low-pass filter,
defined on the complex plane, on the coefficients $a_{k}$. We start with
the Taylor coefficients of $w(z)$, which can be written in terms of its
real an imaginary parts,

\begin{displaymath}
  w(z)
  =
  f_{\rm c}(\rho,\theta)
  +
  \ii
  f_{\rm s}(\rho,\theta),
\end{displaymath}

\noindent
in terms of which the Taylor coefficients are given by

\noindent
\begin{eqnarray*}
  a_{k}
  & = &
  \frac{\rho^{-k}}{\pi}
  \int_{-\pi}^{\pi}d\theta\,
  f_{\rm c}(\rho,\theta)
  \cos(k\theta)
  \\
  & = &
  \frac{\rho^{-k}}{\pi}
  \int_{-\pi}^{\pi}d\theta\,
  f_{\rm s}(\rho,\theta)
  \sin(k\theta).
\end{eqnarray*}

\noindent
The Taylor coefficients of $w_{\epsilon}(z)$ are similarly given by

\noindent
\begin{eqnarray*}
  a_{\epsilon,k}
  & = &
  \frac{\rho^{-k}}{\pi}
  \int_{-\pi}^{\pi}d\theta\,
  f_{\epsilon,{\rm c}}(\rho,\theta)
  \cos(k\theta)
  \\
  & = &
  \frac{\rho^{-k}}{\pi}
  \int_{-\pi}^{\pi}d\theta\,
  f_{\epsilon,{\rm s}}(\rho,\theta)
  \sin(k\theta).
\end{eqnarray*}

\noindent
Let us work out only the first case, involving the cosine, since the work
for the second one in essentially identical and leads to the same result.
Using the definition of $f_{\epsilon,{\rm c}}(\rho,\theta)$ in terms of
$f_{\rm c}(\rho,\theta)$ we have

\noindent
\begin{eqnarray*}
  a_{\epsilon,k}
  & = &
  \frac{\rho^{-k}}{\pi}
  \int_{-\pi}^{\pi}d\theta\,
  \cos(k\theta)\,
  \frac{1}{2\epsilon}
  \int_{\theta-\epsilon}^{\theta+\epsilon}d\theta'\,
  f_{\rm c}\!\left(\rho,\theta'\right)
  \\
  & = &
  -\,
  \frac{\rho^{-k}}{2\epsilon \pi}
  \int_{-\pi}^{\pi}d\theta\,
  \frac{\sin(k\theta)}{k}\,
  \frac{d}{d\theta}
  \int_{\theta-\epsilon}^{\theta+\epsilon}d\theta'\,
  f_{\rm c}\!\left(\rho,\theta'\right)
  \\
  & = &
  -\,
  \frac{\rho^{-k}}{2\epsilon \pi k}
  \int_{-\pi}^{\pi}d\theta\,
  \sin(k\theta)\;
  f_{\rm c}\!\left(\rho,\theta'\right)
  \at{\theta-\epsilon}{\theta+\epsilon}
  \\
  & = &
  -\,
  \frac{\rho^{-k}}{2\epsilon \pi k}
  \int_{-\pi}^{\pi}d\theta\,
  \sin(k\theta)
  f_{\rm c}(\rho,\theta+\epsilon)
  +
  \frac{\rho^{-k}}{2\epsilon \pi k}
  \int_{-\pi}^{\pi}d\theta\,
  \sin(k\theta)
  f_{\rm c}(\rho,\theta-\epsilon),
\end{eqnarray*}

\noindent
where we integrated by parts and where there is no integrated term due to
the periodicity of the integrand in $\theta$. We now change variables in
each integral, using $\theta'=\theta\pm\epsilon$, in order to obtain

\noindent
\begin{eqnarray*}
  a_{\epsilon,k}
  & = &
  -\,
  \frac{\rho^{-k}}{2\epsilon \pi k}
  \int_{-\pi}^{\pi}d\theta'\,
  \sin\!\left(k\theta'-k\epsilon\right)
  f_{\rm c}\!\left(\rho,\theta'\right)
  +
  \frac{\rho^{-k}}{2\epsilon \pi k}
  \int_{-\pi}^{\pi}d\theta'\,
  \sin\!\left(k\theta'+k\epsilon\right)
  f_{\rm c}\!\left(\rho,\theta'\right)
  \\
  & = &
  \frac{\rho^{-k}}{2\epsilon \pi k}
  \int_{-\pi}^{\pi}d\theta'\,
  f_{\rm c}\!\left(\rho,\theta'\right)
  \left[
    \sin\!\left(k\theta'+k\epsilon\right)
    -
    \sin\!\left(k\theta'-k\epsilon\right)
  \right],
\end{eqnarray*}

\noindent
where the integration limits did not change in the transformations of
variables due to the periodicity of the integrand. Changing $\theta'$ back
to $\theta$ we are left with

\noindent
\begin{eqnarray*}
  a_{\epsilon,k}
  & = &
  \frac{\rho^{-k}}{2\epsilon \pi k}
  \int_{-\pi}^{\pi}d\theta\,
  f_{\rm c}\!\left(\rho,\theta\right)
  \left[
    \sin(k\theta)
    \cos(k\epsilon)
    +
    \sin(k\epsilon)
    \cos(k\theta)
  \right.
  +
  \\
  &   &
  \hspace{9em}
  -
  \left.
    \sin(k\theta)
    \cos(k\epsilon)
    +
    \sin(k\epsilon)
    \cos(k\theta)
  \right]
  \\
  & = &
  \left[
    \frac{\sin(k\epsilon)}{(k\epsilon)}
  \right]
  \frac{\rho^{-k}}{\pi}
  \int_{-\pi}^{\pi}d\theta\,
  f_{\rm c}\!\left(\rho,\theta\right)
  \cos(k\theta).
\end{eqnarray*}

\noindent
Since we recover in this way the expression of the coefficients of $f_{\rm
  c}(\rho,\theta)$, we get

\begin{displaymath}
  a_{\epsilon,k}
  =
  \left[
    \frac{\sin(k\epsilon)}{(k\epsilon)}
  \right]
  a_{k},
\end{displaymath}

\noindent
which is the same result obtained in the text through the application of
the filter, as an operator, to the expansion of $w(z)$. This more direct
derivation bypasses any preoccupations with the convergence of the series
during that process, due to the term-wise application of the integral
operator.

\subsection{Alternate Proof of the Inner Analyticity of {\boldmath
    $w_{\epsilon}(z)$}}

Here we establish that $w_{\epsilon}(z)$ is analytic by showing that its
real and imaginary parts satisfy the Cauchy-Riemann conditions. Consider
an inner analytic function $w(z)$ and the corresponding filtered function
within the open unit disk, with the real angular parameter
$0<\epsilon\leq\pi$

\begin{displaymath}
  w_{\epsilon}(z)
  =
  \frac{1}{2\epsilon}
  \int_{\theta-\epsilon}^{\theta+\epsilon}d\theta'\,
  w\!\left(z'\right).
\end{displaymath}

\noindent
Since $w(z)$ is analytic, we have $w(z)=f_{\rm c}(\rho,\theta)+\ii f_{\rm
  s}(\rho,\theta)$ where $f_{\rm c}(\rho,\theta)$ and $f_{\rm
  s}(\rho,\theta)$ satisfy the Cauchy-Riemann conditions in polar
coordinates,

\noindent
\begin{eqnarray*}
  \frac{\partial f_{\rm c}}{\partial\rho}(\rho,\theta)
  & = &
  \frac{1}{\rho}\,
  \frac{\partial f_{\rm s}}{\partial\theta}(\rho,\theta),
  \\
  \frac{1}{\rho}\,
  \frac{\partial f_{\rm c}}{\partial\theta}(\rho,\theta)
  & = &
  -\,
  \frac{\partial f_{\rm s}}{\partial\rho}(\rho,\theta).
\end{eqnarray*}

\noindent
It follows that the filtered function can be written as

\noindent
\begin{eqnarray*}
  w_{\epsilon}(z)
  & = &
  f_{\epsilon,{\rm c}}(\rho,\theta)+\ii f_{\epsilon,{\rm s}}(\rho,\theta)
  \\
  & = &
  \frac{1}{2\epsilon}
  \int_{\theta-\epsilon}^{\theta+\epsilon}d\theta'\,
  f_{\rm c}\!\left(\rho,\theta'\right)
  +
  \ii\,
  \frac{1}{2\epsilon}
  \int_{\theta-\epsilon}^{\theta+\epsilon}d\theta'\,
  f_{\rm s}\!\left(\rho,\theta'\right),
\end{eqnarray*}

\noindent
so that we have

\noindent
\begin{eqnarray*}
  f_{\epsilon,{\rm c}}(\rho,\theta)
  & = &
  \frac{1}{2\epsilon}
  \int_{\theta-\epsilon}^{\theta+\epsilon}d\theta'\,
  f_{\rm c}\!\left(\rho,\theta'\right),
  \\
  f_{\epsilon,{\rm s}}(\rho,\theta)
  & = &
  \frac{1}{2\epsilon}
  \int_{\theta-\epsilon}^{\theta+\epsilon}d\theta'\,
  f_{\rm s}\!\left(\rho,\theta'\right).
\end{eqnarray*}

\noindent
Since $f_{\rm c}(\rho,\theta)$ and $f_{\rm s}(\rho,\theta)$ are continuous
and differentiable, it is clear that so are $f_{\epsilon,{\rm
    c}}(\rho,\theta)$ and $f_{\epsilon,{\rm s}}(\rho,\theta)$. If we
calculate their partial derivatives with respect to $\rho$ we get

\noindent
\begin{eqnarray*}
  \frac{\partial f_{\epsilon,{\rm c}}}{\partial\rho}(\rho,\theta)
  & = &
  \frac{1}{2\epsilon}
  \int_{\theta-\epsilon}^{\theta+\epsilon}d\theta'\,
  \frac{\partial f_{\rm c}}{\partial\rho}\!\left(\rho,\theta'\right),
  \\
  \frac{\partial f_{\epsilon,{\rm s}}}{\partial\rho}(\rho,\theta)
  & = &
  \frac{1}{2\epsilon}
  \int_{\theta-\epsilon}^{\theta+\epsilon}d\theta'\,
  \frac{\partial f_{\rm s}}{\partial\rho}\!\left(\rho,\theta'\right).
\end{eqnarray*}

\noindent
Using the Cauchy-Riemann relations for $w(z)$ we may write these as

\noindent
\begin{eqnarray*}
  \frac{\partial f_{\epsilon,{\rm c}}}{\partial\rho}(\rho,\theta)
  & = &
  \frac{1}{2\epsilon}
  \int_{\theta-\epsilon}^{\theta+\epsilon}d\theta'\,
  \frac{1}{\rho}\,
  \frac{\partial f_{\rm s}}{\partial\theta}\!\left(\rho,\theta'\right)
  \\
  & = &
  \frac{1}{\rho}\,
  \frac{1}{2\epsilon}\,
  f_{\rm s}\!\left(\rho,\theta'\right)
  \at{\theta-\epsilon}{\theta+\epsilon},
  \\
  & = &
  \frac{1}{\rho}\,
  \frac
  {f_{\rm s}(\rho,\theta+\epsilon)-f_{\rm s}(\rho,\theta-\epsilon)}
  {2\epsilon},
  \\
  \frac{\partial f_{\epsilon,{\rm s}}}{\partial\rho}(\rho,\theta)
  & = &
  -\,
  \frac{1}{2\epsilon}
  \int_{\theta-\epsilon}^{\theta+\epsilon}d\theta'\,
  \frac{1}{\rho}\,
  \frac{\partial f_{\rm c}}{\partial\theta}\!\left(\rho,\theta'\right)
  \\
  & = &
  -\,
  \frac{1}{\rho}\,
  \frac{1}{2\epsilon}\,
  f_{\rm c}\!\left(\rho,\theta'\right)
  \at{\theta-\epsilon}{\theta+\epsilon}
  \\
  & = &
  -\,
  \frac{1}{\rho}\,
  \frac
  {f_{\rm c}(\rho,\theta+\epsilon)-f_{\rm c}(\rho,\theta-\epsilon)}
  {2\epsilon}.
\end{eqnarray*}

\noindent
If we now calculate the partial derivatives of $f_{\epsilon,{\rm
    c}}(\rho,\theta)$ and $f_{\epsilon,{\rm s}}(\rho,\theta)$ with respect
to $\theta$ we get

\noindent
\begin{eqnarray*}
  \frac{1}{\rho}\,
  \frac{\partial f_{\epsilon,{\rm c}}}{\partial\theta}(\rho,\theta)
  & = &
  \frac{1}{\rho}\,
  \frac{1}{2\epsilon}\,
  \frac{\partial }{\partial\theta}
  \int_{\theta-\epsilon}^{\theta+\epsilon}d\theta'\,
  f_{\rm c}\!\left(\rho,\theta'\right)
  \\
  & = &
  \frac{1}{\rho}\,
  \frac
  {f_{\rm c}(\rho,\theta+\epsilon)-f_{\rm c}(\rho,\theta-\epsilon)}
  {2\epsilon},
  \\
  \frac{1}{\rho}\,
  \frac{\partial f_{\epsilon,{\rm s}}}{\partial\theta}(\rho,\theta)
  & = &
  \frac{1}{\rho}\,
  \frac{1}{2\epsilon}\,
  \frac{\partial }{\partial\theta}
  \int_{\theta-\epsilon}^{\theta+\epsilon}d\theta'\,
  f_{\rm s}\!\left(\rho,\theta'\right)
  \\
  & = &
  \frac{1}{\rho}\,
  \frac
  {f_{\rm s}(\rho,\theta+\epsilon)-f_{\rm s}(\rho,\theta-\epsilon)}
  {2\epsilon}.
\end{eqnarray*}

\noindent
Comparing this pair of equations with the previous one we get

\noindent
\begin{eqnarray*}
  \frac{\partial f_{\epsilon,{\rm c}}}{\partial\rho}(\rho,\theta)
  & = &
  \frac{1}{\rho}\,
  \frac{\partial f_{\epsilon,{\rm s}}}{\partial\theta}(\rho,\theta),
  \\
  \frac{1}{\rho}\,
  \frac{\partial f_{\epsilon,{\rm c}}}{\partial\theta}(\rho,\theta)
  & = &
  -\,
  \frac{\partial f_{\epsilon,{\rm s}}}{\partial\rho}(\rho,\theta),
\end{eqnarray*}

\noindent
which establish the analyticity of $w_{\epsilon}(z)$, in the same domain
as that of $w(z)$. Let us now examine the other properties defining an
inner analytic function. For one thing we have $w(0)=0$, which means that

\noindent
\begin{eqnarray*}
  \lim_{\rho\to 0}
  f_{\rm c}(\rho,\theta)
  & = &
  0,
  \\
  \lim_{\rho\to 0}
  f_{\rm s}(\rho,\theta)
  & = &
  0.
\end{eqnarray*}

\noindent
If we calculate the corresponding limits for $w_{\epsilon}(z)$ we get

\noindent
\begin{eqnarray*}
  \lim_{\rho\to 0}
  f_{\epsilon,{\rm c}}(\rho,\theta)
  & = &
  \lim_{\rho\to 0}
  \frac{1}{2\epsilon}
  \int_{\theta-\epsilon}^{\theta+\epsilon}d\theta'\,
  f_{\rm c}\!\left(\rho,\theta'\right)
  \\
  & = &
  \frac{1}{2\epsilon}
  \int_{\theta-\epsilon}^{\theta+\epsilon}d\theta'\,
  \lim_{\rho\to 0}
  f_{\rm c}\!\left(\rho,\theta'\right)
  \\
  & = &
  0,
  \\
  \lim_{\rho\to 0}
  f_{\epsilon,{\rm s}}(\rho,\theta)
  & = &
  \lim_{\rho\to 0}
  \frac{1}{2\epsilon}
  \int_{\theta-\epsilon}^{\theta+\epsilon}d\theta'\,
  f_{\rm s}\!\left(\rho,\theta'\right)
  \\
  & = &
  \frac{1}{2\epsilon}
  \int_{\theta-\epsilon}^{\theta+\epsilon}d\theta'\,
  \lim_{\rho\to 0}
  f_{\rm s}\!\left(\rho,\theta'\right)
  \\
  & = &
  0.
\end{eqnarray*}

\noindent
We have therefore that $w_{\epsilon}(0)=0$. Finally, $w(z)$ reduces to a
real function over the interval $(-1,1)$ of the real axis, which means
that its imaginary part is zero there, and therefore that $f_{\rm
  s}(\rho,0)=0$ and $f_{\rm s}(\rho,\pm\pi)=0$. If we write
$f_{\epsilon,{\rm s}}(\rho,\theta)$ for the same values of $\theta$ we get

\noindent
\begin{eqnarray*}
  f_{\epsilon,{\rm s}}(\rho,0)
  & = &
  \frac{1}{2\epsilon}
  \int_{-\epsilon}^{\epsilon}d\theta'\,
  f_{\rm s}\!\left(\rho,\theta'\right),
  \\
  f_{\epsilon,{\rm s}}(\rho,\pm\pi)
  & = &
  \frac{1}{2\epsilon}
  \int_{\pm\pi-\epsilon}^{\pm\pi+\epsilon}d\theta'\,
  f_{\rm s}\!\left(\rho,\theta'\right)
  \\
  & = &
  \frac{1}{2\epsilon}
  \int_{-\pi-\epsilon}^{\pi+\epsilon}d\theta'\,
  f_{\rm s}\!\left(\rho,\theta'\right).
\end{eqnarray*}

\noindent
However, since $w(z)$ is an inner analytic function we have that $f_{\rm
  s}(\rho,\theta)$ is an odd function of $\theta$. In both cases above
this implies that the integral is zero, and hence we conclude that
$w_{\epsilon}(z)$ is an inner analytic function as well.

\subsection{Proof of Analyticity of {\boldmath $w_{\varepsilon}(z)$} in
  the Cartesian Case}

As a curiosity, it is interesting to point out that a first-order linear
low-pass filter can be defined over a straight segment on the complex
plane. In this way a filter in the Cartesian coordinates of the complex
plane can be defined. Consider an analytic function $w(z)$ anywhere on the
complex plane. Consider also a segment of length $2\varepsilon$ and a
fixed direction given by a constant angle $\alpha$ with the real
axis. Given an arbitrary position $z$ on the complex plane we define then
two other points by

\noindent
\begin{eqnarray*}
  z_{\oplus}
  & = &
  z+\varepsilon\e{\ii\alpha},
  \\
  z_{\ominus}
  & = &
  z-\varepsilon\e{\ii\alpha}.
\end{eqnarray*}

\noindent
This defines a segment of length $2\varepsilon$ going from $z_{\ominus}$
to $z_{\oplus}$. Given any point $z$ such that this segment is contained
within the analyticity domain of $w(z)$, we may now define a filtered
function $w_{\varepsilon}(z)$ by

\noindent
\begin{eqnarray*}
  w_{\varepsilon}(z)
  & = &
  \frac{1}{2\varepsilon}
  \int_{z_{\ominus}}^{z_{\oplus}}dz'\,
  w\!\left(z'\right)
  \\
  & = &
  \frac{\e{\ii\alpha}}{2\varepsilon}
  \int_{-\varepsilon}^{\varepsilon}d\lambda\,
  w\!\left(z+\lambda\e{\ii\alpha}\right),
\end{eqnarray*}

\noindent
where $\lambda$ is a real parameter describing the segment, such that
$-\varepsilon\leq\lambda\leq\varepsilon$ and

\noindent
\begin{eqnarray*}
  z'
  & = &
  z+\lambda\e{\ii\alpha}
  \;\;\;\Rightarrow
  \\
  dz'
  & = &
  \e{\ii\alpha}d\lambda.
\end{eqnarray*}

\noindent
Since $w(z)=u(x,y)+\ii v(x,y)$ is analytic, $u(x,y)$ and $v(x,y)$ are
continuous, differentiable and satisfy the Cauchy-Riemann conditions in
Cartesian coordinates. We may write for $w_{\varepsilon}(z)$

\noindent
\begin{eqnarray*}
  w_{\varepsilon}(z)
  & = &
  u_{\varepsilon}(x,y)
  +
  \ii
  v_{\varepsilon}(x,y)
  \\
  & = &
  \frac{\e{\ii\alpha}}{2\varepsilon}
  \int_{-\varepsilon}^{\varepsilon}d\lambda\,
  \left\{
    \rule{0em}{2.5ex}
    u[x+\lambda\cos(\alpha),y+\lambda\sin(\alpha)]
  \right.
  +
  \\
  &   &
  \hspace{5em}
  +
  \left.
    \ii
    v[x+\lambda\cos(\alpha),y+\lambda\sin(\alpha)]
    \rule{0em}{2.5ex}
  \right\}
  \;\;\;\Rightarrow
  \\
  u_{\varepsilon}(x,y)
  & = &
  \frac{\e{\ii\alpha}}{2\varepsilon}
  \int_{-\varepsilon}^{\varepsilon}d\lambda\,
  u[x+\lambda\cos(\alpha),y+\lambda\sin(\alpha)],
  \\
  v_{\varepsilon}(x,y)
  & = &
  \frac{\e{\ii\alpha}}{2\varepsilon}
  \int_{-\varepsilon}^{\varepsilon}d\lambda\,
  v[x+\lambda\cos(\alpha),y+\lambda\sin(\alpha)].
\end{eqnarray*}

\noindent
It is now clear that $u_{\varepsilon}(x,y)$ and $v_{\varepsilon}(x,y)$ are
also continuous and differentiable. If we now take the partial derivatives
of these functions with respect to $x$ we get

\noindent
\begin{eqnarray*}
  \frac{\partial u_{\varepsilon}}{\partial x}(x,y)
  & = &
  \frac{\e{\ii\alpha}}{2\varepsilon}
  \int_{-\varepsilon}^{\varepsilon}d\lambda\,
  \frac{\partial u}{\partial x}[x+\lambda\cos(\alpha),y+\lambda\sin(\alpha)],
  \\
  \frac{\partial v_{\varepsilon}}{\partial x}(x,y)
  & = &
  \frac{\e{\ii\alpha}}{2\varepsilon}
  \int_{-\varepsilon}^{\varepsilon}d\lambda\,
  \frac{\partial v}{\partial x}[x+\lambda\cos(\alpha),y+\lambda\sin(\alpha)].
\end{eqnarray*}

\noindent
Using the Cauchy-Riemann conditions for $u(x,y)$ and $v(x,y)$ we get

\noindent
\begin{eqnarray*}
  \frac{\partial u_{\varepsilon}}{\partial x}(x,y)
  & = &
  \frac{\e{\ii\alpha}}{2\varepsilon}
  \int_{-\varepsilon}^{\varepsilon}d\lambda\,
  \frac{\partial v}{\partial y}[x+\lambda\cos(\alpha),y+\lambda\sin(\alpha)],
  \\
  \frac{\partial v_{\varepsilon}}{\partial x}(x,y)
  & = &
  -\,
  \frac{\e{\ii\alpha}}{2\varepsilon}
  \int_{-\varepsilon}^{\varepsilon}d\lambda\,
  \frac{\partial u}{\partial y}[x+\lambda\cos(\alpha),y+\lambda\sin(\alpha)].
\end{eqnarray*}

\noindent
Taking now the partial derivatives of $u_{\varepsilon}(x,y)$ and
$v_{\varepsilon}(x,y)$ with respect to $y$ we get

\noindent
\begin{eqnarray*}
  \frac{\partial u_{\varepsilon}}{\partial y}(x,y)
  & = &
  \frac{\e{\ii\alpha}}{2\varepsilon}
  \int_{-\varepsilon}^{\varepsilon}d\lambda\,
  \frac{\partial u}{\partial y}[x+\lambda\cos(\alpha),y+\lambda\sin(\alpha)],
  \\
  \frac{\partial v_{\varepsilon}}{\partial y}(x,y)
  & = &
  \frac{\e{\ii\alpha}}{2\varepsilon}
  \int_{-\varepsilon}^{\varepsilon}d\lambda\,
  \frac{\partial v}{\partial y}[x+\lambda\cos(\alpha),y+\lambda\sin(\alpha)].
\end{eqnarray*}

\noindent
Comparing this pair of equation with the previous ones we finally get

\noindent
\begin{eqnarray*}
  \frac{\partial u_{\varepsilon}}{\partial x}(x,y)
  & = &
  \frac{\partial v_{\varepsilon}}{\partial y}(x,y),
  \\
  \frac{\partial u_{\varepsilon}}{\partial y}(x,y)
  & = &
  -\,
  \frac{\partial v_{\varepsilon}}{\partial x}(x,y).
\end{eqnarray*}

\noindent
This establishes that $u_{\varepsilon}(x,y)$ and $v_{\varepsilon}(x,y)$
satisfy the Cauchy-Riemann conditions, and therefore that
$w_{\varepsilon}(z)$ is analytic. Once $w_{\varepsilon}(z)$ is defined by
the filter at all points of the domain of analyticity of $w(z)$ where the
segment fits, and now that it has been proven analytic there, one can
extend the definition of $w_{\varepsilon}(z)$ to the whole domain of
analyticity of $w(z)$ by analytic continuation.

\section{Proof of Convergence to the Infinite-Order
  Kernel}\label{APPInfOrd}

In this appendix we will offer proof of the convergence of the sequence of
order-$N$ scaled kernels to the infinite-order scaled kernel, in the
$N\to\infty$ limit. The point is to show that the infinite sequence of
real functions $\bar{K}_{\epsilon_{N}}^{(N)}(\theta)$, with
$\epsilon<\pi$, converges in the $N\to\infty$ limit to a definite regular
real function with finite support within $[-\epsilon,\epsilon]$, which we
denote as $\bar{K}_{\epsilon}^{(\infty)}(\theta)$. In addition to this, we
will establish that the infinite-order scaled kernel
$\bar{K}_{\epsilon}^{(\infty)}(\theta)$ is a $C^{\infty}$ function, and
that it is not an analytic function in the real sense of the term.

In order to do this we must first establish a few more properties of the
first-order filter, in addition to those demonstrated in~\cite{LPFFSaPDE}.
We will also have to examine more closely some aspects of the process of
construction of the infinite-order scaled kernel.

\subsection{Invariance of the Sign of the First Derivative}

According to one of the properties established before for the first-order
filter~\cite{LPFprop03}, if $f(\theta)$ is continuous in the domain where
the filter is applied, then $f_{\epsilon}(\theta)$ is differentiable, and
its derivative is given by

\begin{equation}\label{firstderiv}
  \frac{df_{\epsilon}(\theta)}{d\theta}
  =
  \frac
  { 
    f\!\left(\theta+\epsilon\right)
    -
    f\!\left(\theta-\epsilon\right)
  }
  {2\epsilon}.
\end{equation}

\noindent
This is valid so long as the support interval of the filter fits
completely inside the region where $f(\theta)$ is continuous. This
immediately implies that, in a region where $f(\theta)$ increases
monotonically we have

\noindent
\begin{eqnarray*}
  f\!\left(\theta+\epsilon\right)
  & \geq &
  f\!\left(\theta-\epsilon\right)
  \;\;\;\Rightarrow
  \\
  \frac{df_{\epsilon}(\theta)}{d\theta}
  & \geq &
  0.
\end{eqnarray*}

\noindent
We may therefore conclude that $f_{\epsilon}(\theta)$ also increases
monotonically within the sub-region where the support interval of the
filter fits inside the region in which $f(\theta)$ is continuous. In the
same way, in a region where $f(\theta)$ decreases monotonically we have

\noindent
\begin{eqnarray*}
  f\!\left(\theta+\epsilon\right)
  & \leq &
  f\!\left(\theta-\epsilon\right)
  \;\;\;\Rightarrow
  \\
  \frac{df_{\epsilon}(\theta)}{d\theta}
  & \leq &
  0.
\end{eqnarray*}

\noindent
We may therefore conclude that $f_{\epsilon}(\theta)$ also decreases
monotonically within that same sub-region. In other words, the monotonic
character of the variation of a function is invariant by the action of the
filter. In particular, at points where $f(\theta)$ is differentiable the
sign of its derivative is invariant by the action of the filter.

\subsection{Invariance of the Sign of the Second Derivative}

If we assume that $f(\theta)$ is differentiable in the domain where the
filter is applied, then $f_{\epsilon}(\theta)$ can be differentiated
twice, and we may obtain its second derivative by simply differentiating
once Equation~(\ref{firstderiv}), which results in

\begin{displaymath}
  \frac{d^{2}f_{\epsilon}(\theta)}{d\theta^{2}}
  =
  \frac{1}{2\epsilon}
  \left[
    \frac{df}{d\theta}\!\left(\theta+\epsilon\right)
    -
    \frac{df}{d\theta}\!\left(\theta-\epsilon\right)
  \right].
\end{displaymath}

\noindent
This is valid so long as the support interval of the filter fits
completely inside the region where $f(\theta)$ is differentiable. This
immediately implies that, in a region where the derivative of $f(\theta)$
increases monotonically we have

\noindent
\begin{eqnarray*}
  \frac{df}{d\theta}\!\left(\theta+\epsilon\right)
  & \geq &
  \frac{df}{d\theta}\!\left(\theta-\epsilon\right)
  \;\;\;\Rightarrow
  \\
  \frac{d^{2}f_{\epsilon}(\theta)}{d\theta^{2}}
  & \geq &
  0.
\end{eqnarray*}

\noindent
We may therefore conclude that the derivative of $f_{\epsilon}(\theta)$
also increases monotonically within the sub-region where the support
interval of the filter fits inside the region in which $f(\theta)$ is
differentiable. In the same way, in a region where the derivative of
$f(\theta)$ decreases monotonically we have

\noindent
\begin{eqnarray*}
  \frac{df}{d\theta}\!\left(\theta+\epsilon\right)
  & \leq &
  \frac{df}{d\theta}\!\left(\theta-\epsilon\right)
  \;\;\;\Rightarrow
  \\
  \frac{d^{2}f_{\epsilon}(\theta)}{d\theta^{2}}
  & \leq &
  0.
\end{eqnarray*}

\noindent
We may therefore conclude that the derivative of $f_{\epsilon}(\theta)$
also decreases monotonically within that same sub-region. In other words,
the monotonic character of the variation of the derivative of a function
is invariant by the action of the filter. In particular, at points where
$f(\theta)$ is twice differentiable the sign of its second derivative is
invariant by the action of the filter.

This implies that in regions where the second derivative of $f(\theta)$
has constant sign, and therefore the concavity of its graph is turned in a
definite direction, up or down, the action of the filter keeps that
concavity turned in the same direction. In other words, away from
inflection points, in regions where the graph of $f(\theta)$ has a
definite concavity turned in a definite direction, $f_{\epsilon}(\theta)$
has the same concavity, turned in the same direction.

\subsection{Action of the Filter in Regions of Definite Concavity}

Consider the action of the first-order filter in a region where the graph
of $f(\theta)$ has definite concavity, turned in a definite direction, and
within which the support interval of the filter fits. As shown
in~\cite{LPFprop01}, if $f(\theta)$ happens to be a linear function within
the support of the filter around a given point, then the filter is the
identity and therefore $f_{\epsilon}(\theta)=f(\theta)$ at that point. On
the other hand, if $f(\theta)$ is not a linear function and its concavity
is turned down, them some values of the function within the support of the
filter must be smaller that in the case of the linear function. Since the
filtered function $f_{\epsilon}(\theta)$ is defined as an average of the
values of $f(\theta)$, it follows that, if the concavity of the graph of
$f(\theta)$ is turned down, then

\begin{displaymath}
  f_{\epsilon}(\theta)
  <
  f(\theta).
\end{displaymath}

\noindent
In the same way, we may conclude that if the concavity of the graph of
$f(\theta)$ is turned up, then

\begin{displaymath}
  f_{\epsilon}(\theta)
  >
  f(\theta).
\end{displaymath}

\noindent
In other words, in regions where the function $f(\theta)$ has its
concavity turned in a definite direction, and within which the support
interval of the filter fits, the action of the filter always changes the
value of the function in the direction to which its concavity is turned.

\subsection{Bounds of the Scaled Kernels}

\begin{figure}[ht]
  \centering
  \fbox{
    \epsfig{file=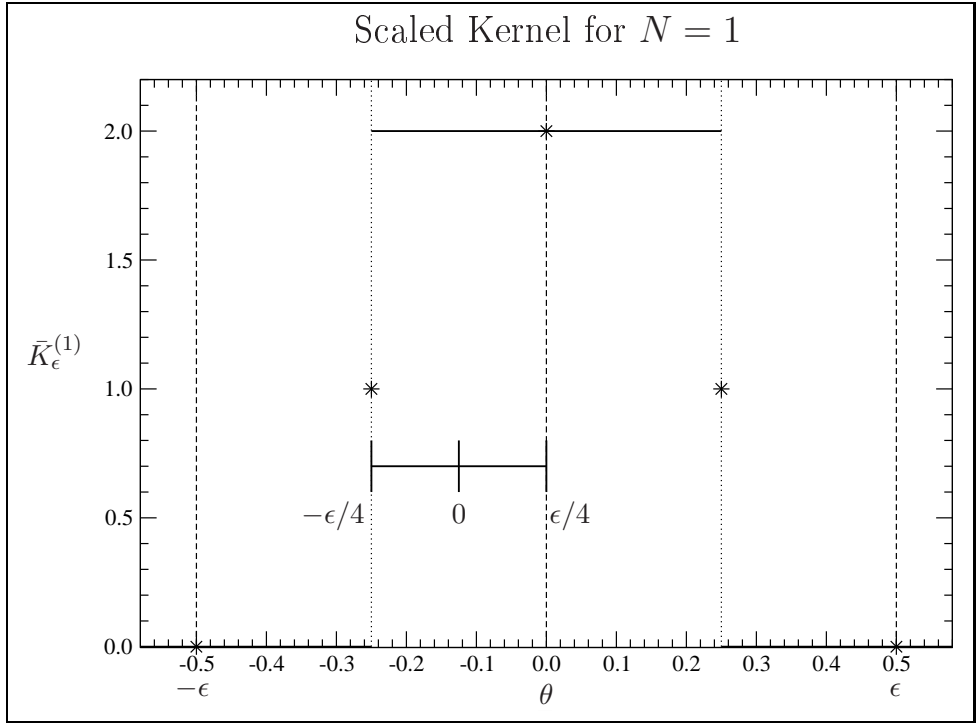,scale=1.0,angle=0}
  }
  \caption{The scaled kernel for $\epsilon=0.5$ and $N=1$, plotted as a
    function of $\theta$ over the support interval $[-\epsilon,\epsilon]$.
    The five invariant points are marked by stars, and the support
    interval of the next filter in the construction sequence is shown. The
    dashed lines mark the center and the two ends of the support interval,
    and hence three of the five invariant points. The dotted lines mark
    the sectors where the function is defined in a piece-wise fashion.}
  \label{Fig12}
\end{figure}

Let us observe that since the filtered function $f_{\epsilon}(\theta)$ is
defined as an average of the function $f(\theta)$, it can never assume
values which are larger than the maximum of the function it is applied on,
or smaller than its minimum, without regard to the value of the range
$\epsilon$. Therefore, since the first kernel we start with in the process
of construction of the infinite-order scaled kernel, that is the kernel
$\bar{K}_{\epsilon_{1}}^{(1)}(\theta)$, with $N=1$ and range $\epsilon/2$,
is bound within the interval $[0,1/\epsilon]$ for all values of $\theta$,
so is the next one, the kernel $\bar{K}_{\epsilon_{2}}^{(2)}(\theta)$. We
may now apply the same argument to this second kernel, and conclude that
the third one in the sequence is also bound in the same way, and so on.
It follows that, for all values of $N$, we have

\begin{displaymath}
  0
  \leq
  \bar{K}_{\epsilon_{N}}^{(N)}(\theta)
  \leq
  \frac{1}{\epsilon},
\end{displaymath}

\noindent
for all values of $\theta$ within the periodic interval $[-\pi,\pi]$, and
in particular for all values of $\theta$ within the support interval
$[-\epsilon,\epsilon]$. It also follows that, if the $N\to\infty$ limit of
the sequence of kernel functions exists, then it is also bound in the same
way.

\subsection{Invariant Points of the Scaled Kernels}

\begin{figure}[ht]
  \centering
  \fbox{
    \epsfig{file=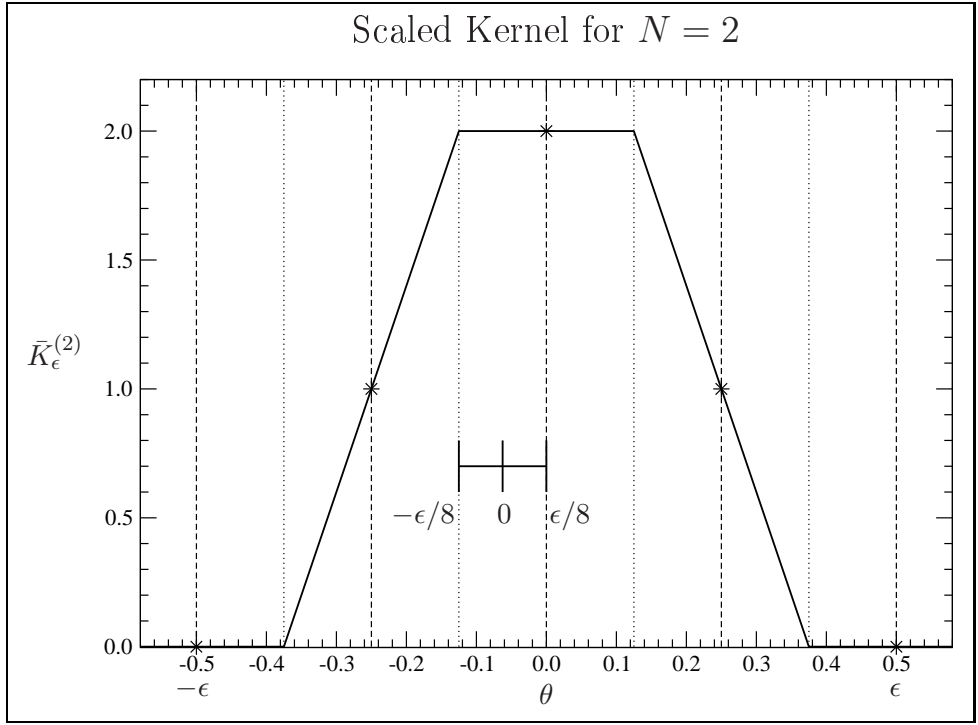,scale=1.0,angle=0}
  }
  \caption{The scaled kernel for $\epsilon=0.5$ and $N=2$, plotted as a
    function of $\theta$ over the support interval $[-\epsilon,\epsilon]$.
    The five invariant points are marked by stars, and the support
    interval of the next filter in the construction sequence is shown. The
    dashed lines mark the four intervals defined by the five invariant
    points. The dotted lines mark the sectors where the function is
    defined in a piece-wise fashion.}
  \label{Fig13}
\end{figure}

\begin{figure}[ht]
  \centering
  \fbox{
    \epsfig{file=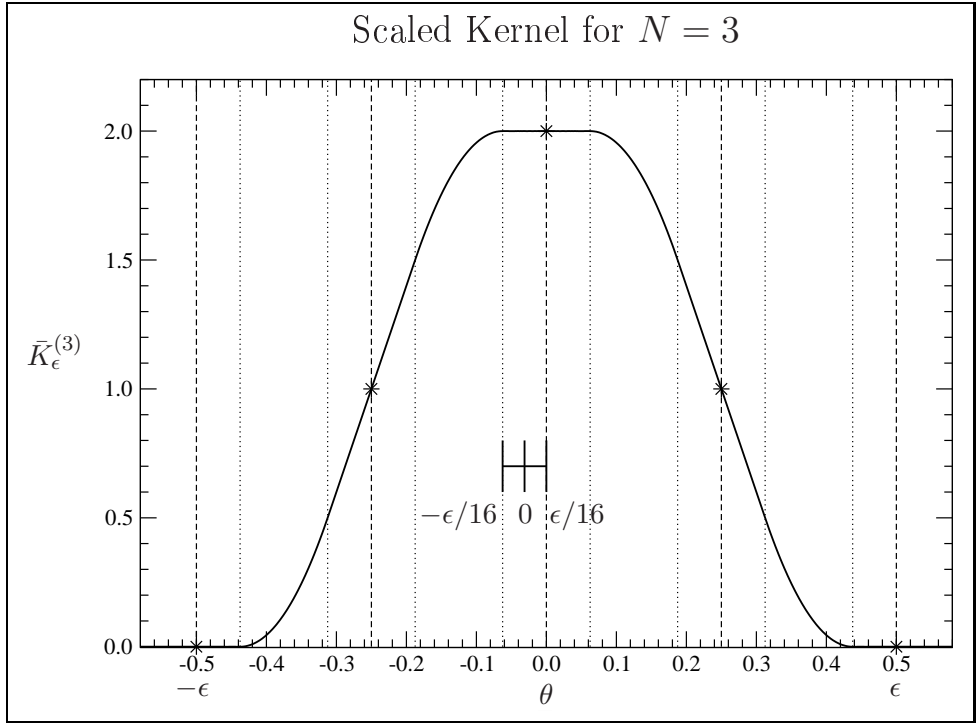,scale=1.0,angle=0}
  }
  \caption{The scaled kernel for $\epsilon=0.5$ and $N=3$, plotted as a
    function of $\theta$ over the support interval $[-\epsilon,\epsilon]$.
    The five invariant points are marked by stars, and the support
    interval of the next filter in the construction sequence is shown. The
    dashed lines mark the four intervals defined by the five invariant
    points. The dotted lines mark the sectors where the function is
    defined in a piece-wise fashion.}
  \label{Fig14}
\end{figure}

Let us now show that there are five points of the graphs of the order-$N$
scaled kernels that remain invariant throughout the construction of the
infinite-order scaled kernel. These are the following: the point of
maximum at $\theta=0$, where the value of the kernel function is
$1/\epsilon$; the two points of minimum at $\theta=\pm\epsilon$, where the
value is zero; and the two inflection points at $\theta=\pm\epsilon/2$,
where the value is $1/(2\epsilon)$. These five points are marked with
stars on the graphs in Figures~\ref{Fig12}, \ref{Fig13} and~\ref{Fig14},
which display the first three kernels in the construction sequence. Note
that in the case of the discontinuous $N=1$ kernel in Figure~\ref{Fig12}
we choose the values at the points of discontinuity according to the
criterion of the average of the lateral limits, which defines these two
future points of inflection in a way that is compatible with this
invariance.

Note now that in Figure~\ref{Fig12} the points of maximum and minimum are
within sectors where the kernel function is linear, in intervals of size
$\epsilon$ (or more, in the case of the points of minimum) around these
points. The support of the next filter to be applied, in the sequence
leading to the construction of the infinite-order scaled kernel, is also
shown in the graph. Since this support has length $\epsilon/2$, it fits
into the intervals where the kernel function is linear, in which case it
acts as the identity, according to one of the properties of the
first-order filter~\cite{LPFprop01}. Therefore, after the application of
this next filter intervals of length $\epsilon/2$ will remain around these
three points, where the next kernel function so generated is still linear.

The result of this operation, which is the $N=2$ scaled kernel, is shown
in Figure~\ref{Fig13}. Note that in this case all the five points listed
before have around them intervals of length $\epsilon/2$ where the kernel
function is linear. Once more the support of the next filter to be
applied, in the sequence leading to the construction of the infinite-order
scaled kernel, is shown in the graph. Since this support has length
$\epsilon/4$, it fits into the intervals where the kernel function is
linear, in which case it acts as the identity, so that after its
application intervals of length $\epsilon/4$ will remain around all these
five points, where the next kernel function generated is still linear. In
particular, this will keep the points invariant, since they are within
sectors where the kernel functions are linear and hence where the
first-order filter acts as the identity.

The result of this last operation, which is the $N=3$ scaled kernel, is
shown in Figure~\ref{Fig14}. Once again all five points are within sectors
of length $\epsilon/4$ where the kernel function is linear. Since the
support of the next filter to be applied has now length $\epsilon/8$, once
more it will keep these points invariant. It is now quite clear that both
the length of the linear sectors and the length of the support of the next
filter will be scaled down exponentially during the process of
construction of the infinite-order kernel, with the support being always
half the length of the intervals, and therefore fitting within them. This
establishes that the five points we listed here are in fact invariant
throughout the construction of the infinite-order scaled kernel. In
particular, it follows that these are the values of the infinite-order
scaled kernel function at these five points, and that the sequence of
order-$N$ scaled kernel functions in fact converges to the infinite-order
scaled kernel function at these five points.

\subsection{Convergence of the Scaled Kernels}\label{APPInfOrdLim}

We are now ready to show that the sequence of order-$N$ scaled kernels
converges to the infinite-order scaled kernel, within the whole support
interval. Of course the convergence is guaranteed outside the support
interval, since all the scaled kernels in the construction sequence are
identically zero there. Starting from the $N=3$ scaled kernel shown in
Figure~\ref{Fig14}, which is an everywhere continuous and differentiable
function, we consider the action on it of the next first-order filter.
Observe that within each one of the four intervals of length $\epsilon/2$
defined by the five invariant points this kernel is monotonic, and also
that its first derivative is monotonic as well. The situation is as
follows: in the first interval $[-\epsilon,-\epsilon/2]$ both the kernel
function and its derivative are monotonically increasing; in the second
interval $[-\epsilon/2,0]$ the kernel function is monotonically
increasing, but its derivative is monotonically decreasing; in the third
interval $[0,\epsilon/2]$, both the kernel function and its derivative are
monotonically decreasing; in the fourth interval $[\epsilon/2,\epsilon]$,
the kernel function is monotonically decreasing, but its derivative is
monotonically increasing.

This means that this kernel function has a definite concavity in each of
the four intervals. Let us now consider the action of the next instance of
the first-order filter, at each point within the final support interval
$[-\epsilon,\epsilon]$. We can say that either one of the five invariant
points is contained within the support of the filter, or none is. If one
of them is contained in the support, then we have already established that
the support of the filter is contained within an interval where the kernel
function is linear, and therefore the filter acts as the identity. In this
case the kernel function is not changed at all. Otherwise, the support is
contained within one of the four intervals where the kernel function has a
definite concavity. In this case the kernel function will be changed, but
its monotonic character, and that of its derivative, will be preserved. In
other words the next kernel will have the same monotonicity and concavity
properties on the same four intervals. Since this argument can then be
iterated, we conclude that all subsequent kernel functions in the
construction sequence have these same monotonicity and concavity
properties, on the same four intervals.

Let us now consider the action of the first-order filter at any subsequent
stage of the construction process. Once again, we have that either one of
the five invariant points is contained within the support of the current
filter, or none is. If one of the points is contained in the support, then
the support of the filter is contained within an interval where the
current kernel function is linear, and therefore the filter acts as the
identity, so that the value of the kernel function is not changed. If none
of the five points is contained within the support, then that support is
contained within one of the four intervals where the current kernel has
the same monotonicity and concavity properties of all the others in the
sequence, starting with $N=3$. This means that at all stages of the
construction process the points of the graph of the current kernel will
always be changed in the same direction within these four intervals, being
therefore always increased in the intervals $[-\epsilon,-\epsilon/2]$ and
$[\epsilon/2,\epsilon]$, and always decreased in the intervals
$[-\epsilon/2,0]$ and $[0,\epsilon/2]$.

What we may conclude from this is that, given any value of $\theta$ within
$[-\epsilon,\epsilon]$, either it is one of the invariant points, at which
all order-$N$ kernel functions have the same values, and therefore where
the sequence of kernel functions converges, or it is a point strictly
within one of the four intervals where the kernel functions have definite
monotonicity and concavity properties. In this case the sequence of values
of the order-$N$ kernel functions at that point form a monotonic real
sequence. Since this is a monotonic sequence of real numbers that is bound
from below by zero and from above by $1/\epsilon$, it follows that the
sequence converges. Since we may therefore state that the point-wise
convergence holds for all points within the final support interval
$[-\epsilon,\epsilon]$, and recalling that outside this interval all
order-$N$ kernels are identically zero, we conclude the the sequence of
order-$N$ scales kernels converges, in the $N\to\infty$ limit, to the
infinite-order scaled kernel, a definite limited real function
$\bar{K}_{\epsilon}^{(\infty)}(\theta)$ with compact support.

\subsection{Differentiability of the Infinite-Order Scaled
  Kernel}\label{APPInfOrdDiff}

The infinite-order scaled kernel with finite range $\epsilon$ has an
interesting property of its own, namely that there is a certain similarity
between the kernel and its derivatives. Every finite-order derivative of
the infinite-order scaled kernel function is made out of a certain number
of rescaled copies of the kernel itself, concatenated together. This is a
consequence of the fact that there is a certain relation between the first
derivative of the order-$N$ scaled kernel and the order-$(N-1)$ scaled
kernel. After this relation is established it can be iterated, resulting
in similar relations for the higher-order derivatives. This property
allows one to establish the existence of the $N\to\infty$ limits of all
the finite-order derivatives, and thus to prove that the infinite-order
scaled kernel is differentiable to all orders.

One can derive the relation between the first derivative of the order-$N$
scaled kernel and the order-$(N-1)$ scaled kernel as follows. If we start
with the Fourier expansion of the order-$N$ scaled kernel, written in the
form

\begin{displaymath}
  \bar{K}_{\epsilon_{N}}^{(N)}(\theta)
  =
  \frac{1}{2\pi}
  +
  \frac{1}{\pi}
  \sum_{k=1}^{\infty}
  \left[
    \,
    \prod_{n=1}^{N}
    \frac
    {\sin\!\left(k\epsilon/2^{n}\right)}
    {\left(k\epsilon/2^{n}\right)}
  \right]
  \cos(k\theta),
\end{displaymath}

\noindent
we may differentiate once term-by-term and thus obtain

\noindent
\begin{eqnarray*}
  \frac{d}{d\theta}
  \bar{K}_{\epsilon_{N}}^{(N)}(\theta)
  & = &
  \frac{1}{\pi}
  \sum_{k=1}^{\infty}
  \left[
    \frac
    {\sin\!\left(k\epsilon/2^{1}\right)}
    {\left(k\epsilon/2^{1}\right)}
  \right]
  \left[
    \,
    \prod_{n=2}^{N}
    \frac
    {\sin\!\left(k\epsilon/2^{n}\right)}
    {\left(k\epsilon/2^{n}\right)}
  \right]
  (-k)
  \sin(k\theta)
  \\
  & = &
  \frac{1}{\pi\epsilon}
  \sum_{k=1}^{\infty}
  \left[
    \,
    \prod_{n=2}^{N}
    \frac
    {\sin\!\left(k\epsilon/2^{n}\right)}
    {\left(k\epsilon/2^{n}\right)}
  \right]
  \left[
    -2
    \sin\!\left(k\epsilon/2\right)
    \sin(k\theta)
  \right],
\end{eqnarray*}

\noindent
noting that for sufficiently large $N$ all the series involved are
absolutely and uniformly convergent. By means of simple trigonometric
identities the product of two sines within brackets can now be written as

\begin{displaymath}
  -2
  \sin\!\left(k\epsilon/2\right)
  \sin(k\theta)
  =
  \cos\!\left[k\left(\theta+\epsilon/2\right)\right]
  -
  \cos\!\left[k\left(\theta-\epsilon/2\right)\right],
\end{displaymath}

\noindent
so that we have for the derivative of the kernel

\noindent
\begin{eqnarray*}
  \frac{d}{d\theta}
  \bar{K}_{\epsilon_{N}}^{(N)}(\theta)
  & = &
  \frac{1}{\pi\epsilon}
  \sum_{k=1}^{\infty}
  \left[
    \,
    \prod_{n=2}^{N}
    \frac
    {\sin\!\left(k\epsilon/2^{n}\right)}
    {\left(k\epsilon/2^{n}\right)}
  \right]
  \cos\!\left[k\left(\theta+\epsilon/2\right)\right]
  +
  \\
  &   &
  -\,
  \frac{1}{\pi\epsilon}
  \sum_{k=1}^{\infty}
  \left[
    \,
    \prod_{n=2}^{N}
    \frac
    {\sin\!\left(k\epsilon/2^{n}\right)}
    {\left(k\epsilon/2^{n}\right)}
  \right]
  \cos\!\left[k\left(\theta-\epsilon/2\right)\right].
\end{eqnarray*}

\begin{figure}[ht]
  \centering
  \fbox{
    \epsfig{file=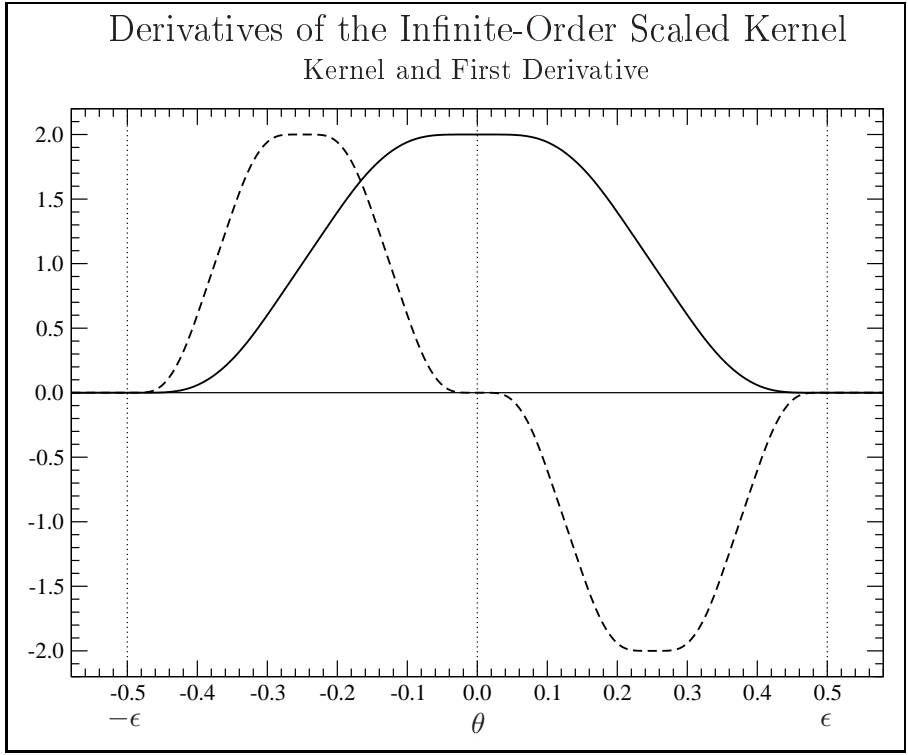,scale=1.0,angle=0}
  }
  \caption{The infinite-order scaled kernel (solid line) and its first
    derivative (dashed line), obtained via the use of their Fourier
    series, for $\epsilon=0.5$, for a large value of $N$ ($100$), plotted
    as functions of $\theta$ over the support interval
    $[-\epsilon,\epsilon]$. The derivative was rescaled down to have the
    same amplitude as the kernel. The dotted lines mark the points where
    the first derivative is zero.}
  \label{Fig15}
\end{figure}

\noindent
If we now define $\epsilon'=\epsilon/2$, we may write

\noindent
\begin{eqnarray*}
  \frac{d}{d\theta}
  \bar{K}_{\epsilon_{N}}^{(N)}(\theta)
  & = &
  \frac{2}{\pi\epsilon'}
  \sum_{k=1}^{\infty}
  \left[
    \,
    \prod_{n=2}^{N}
    \frac
    {\sin\!\left(k\epsilon'/2^{n-1}\right)}
    {\left(k\epsilon'/2^{n-1}\right)}
  \right]
  \cos\!\left[k\left(\theta+\epsilon'\right)\right]
  +
  \\
  &   &
  -\,
  \frac{2}{\pi\epsilon'}
  \sum_{k=1}^{\infty}
  \left[
    \,
    \prod_{n=2}^{N}
    \frac
    {\sin\!\left(k\epsilon'/2^{n-1}\right)}
    {\left(k\epsilon'/2^{n-1}\right)}
  \right]
  \cos\!\left[k\left(\theta-\epsilon'\right)\right]
  \\
  & = &
  \frac{1}{2\pi}
  +
  \frac{2}{\pi\epsilon'}
  \sum_{k=1}^{\infty}
  \left[
    \,
    \prod_{n'=1}^{N-1}
    \frac
    {\sin\!\left(k\epsilon'/2^{n'}\right)}
    {\left(k\epsilon'/2^{n'}\right)}
  \right]
  \cos\!\left[k\left(\theta+\epsilon'\right)\right]
  +
  \\
  &   &
  -\,
  \frac{1}{2\pi}
  -
  \frac{2}{\pi\epsilon'}
  \sum_{k=1}^{\infty}
  \left[
    \,
    \prod_{n'=1}^{N-1}
    \frac
    {\sin\!\left(k\epsilon'/2^{n'}\right)}
    {\left(k\epsilon'/2^{n'}\right)}
  \right]
  \cos\!\left[k\left(\theta-\epsilon'\right)\right],
\end{eqnarray*}

\noindent
where we made $n'=n-1$, which implies $n=n'+1$. We see therefore that in
this way we recover in the right-hand side the expression of the scaled
kernel of order $N-1$ with range $\epsilon/2$, so that we have, writing
$\epsilon'$ back in terms of $\epsilon$,

\begin{displaymath}
  \frac{d}{d\theta}
  \bar{K}_{\epsilon_{N}}^{(N)}(\theta)
  =
  \frac{1}{\epsilon}
  \left[
    \bar{K}_{\epsilon_{(N-1)}/2}^{(N-1)}\!\left(\theta+\epsilon/2\right)
    -
    \bar{K}_{\epsilon_{(N-1)}/2}^{(N-1)}\!\left(\theta-\epsilon/2\right)
  \right].
\end{displaymath}

\begin{figure}[ht]
  \centering
  \fbox{
    \epsfig{file=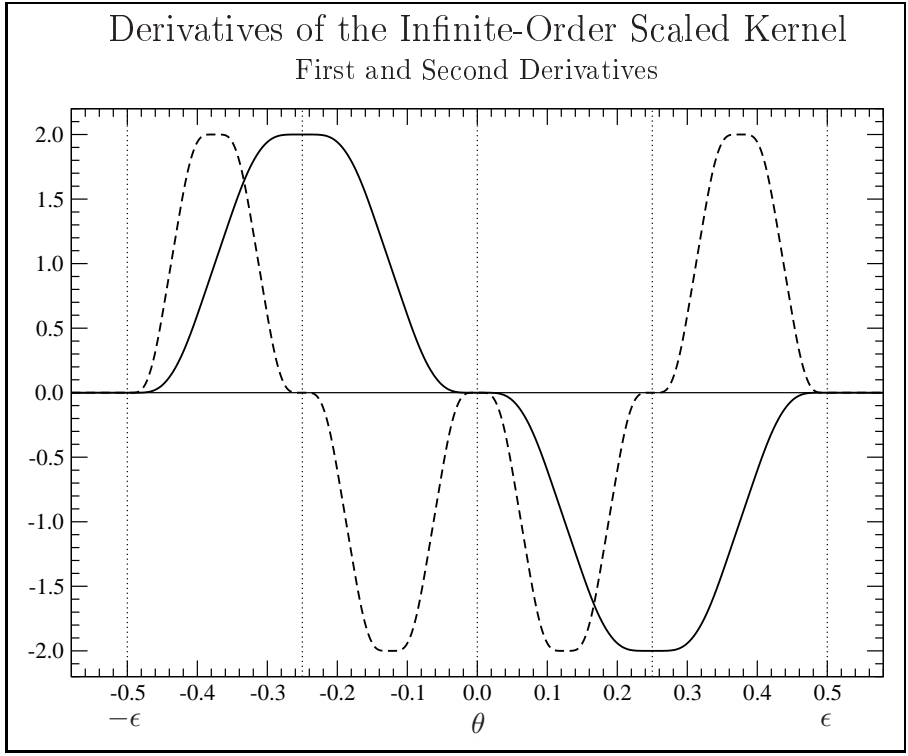,scale=1.0,angle=0}
  }
  \caption{The first (solid line) and second (dashed line) derivatives of
    the infinite-order scaled kernel, obtained via the use of their
    Fourier series, for $\epsilon=0.5$, for a large value of $N$ ($100$),
    plotted as functions of $\theta$ over the support interval
    $[-\epsilon,\epsilon]$. Both derivatives were rescaled down to have
    the same amplitude as the kernel. The dotted lines mark the points
    where the second derivative is zero.}
  \label{Fig16}
\end{figure}

\noindent
Since we already know that the $N\to\infty$ limit of the right-hand side
exists, this establishes that the $N\to\infty$ limit of the left-hand
exists as well. Taking the $N\to\infty$ limit we end up with the
infinite-order scaled kernel on both sides, so that we have the relation

\begin{displaymath}
  \frac{d}{d\theta}
  \bar{K}_{\epsilon}^{\infty}(\theta)
  =
  \frac{1}{\epsilon}
  \left[
    \bar{K}_{\epsilon/2}^{\infty}\!\left(\theta+\epsilon/2\right)
    -
    \bar{K}_{\epsilon/2}^{\infty}\!\left(\theta-\epsilon/2\right)
  \right].
\end{displaymath}

\begin{figure}[ht]
  \centering
  \fbox{
    \epsfig{file=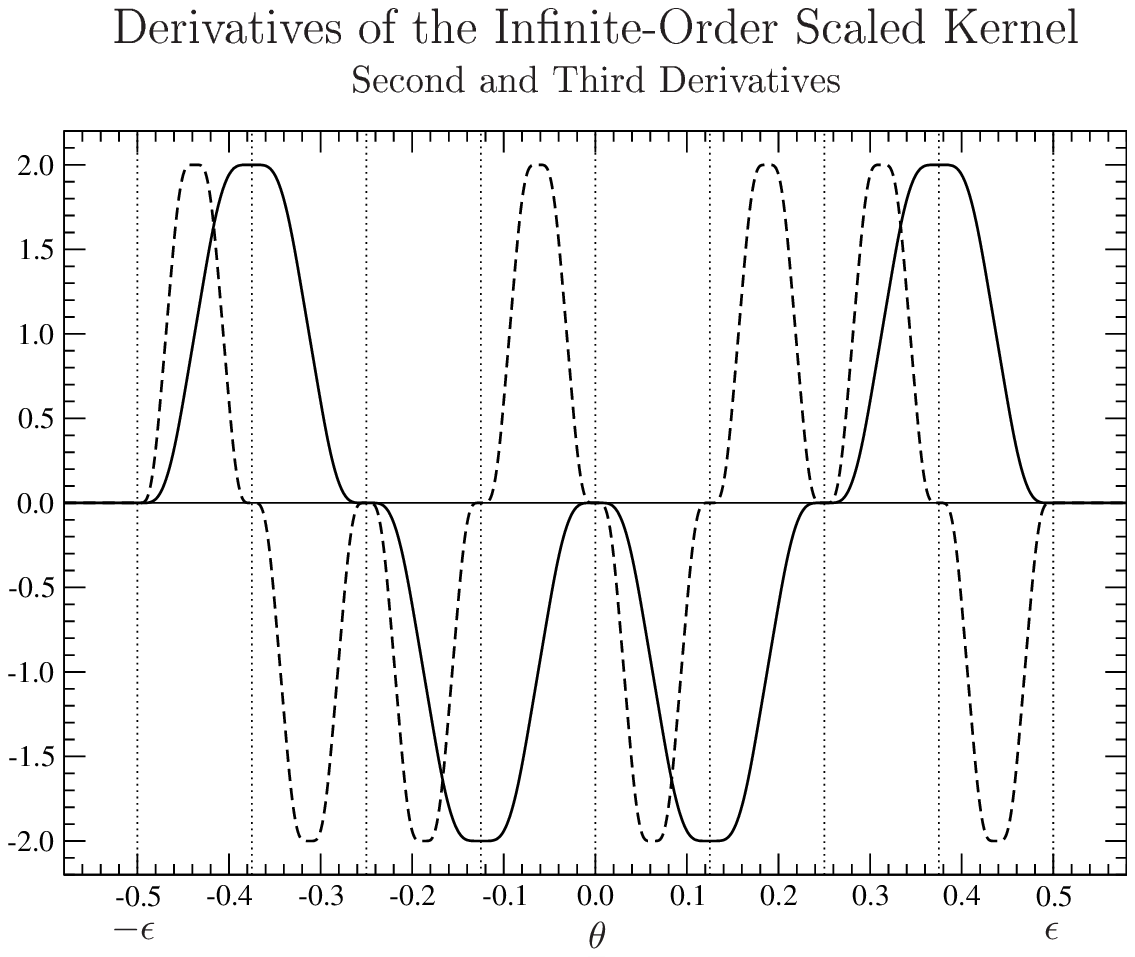,scale=1.0,angle=0}
  }
  \caption{The second (solid line) and third (dashed line) derivatives of
    the infinite-order scaled kernel, obtained via the use of their
    Fourier series, for $\epsilon=0.5$, for a large value of $N$ ($100$),
    plotted as functions of $\theta$ over the support interval
    $[-\epsilon,\epsilon]$. Both derivatives were rescaled down to have
    the same amplitude as the kernel. The dotted lines mark the points
    where the third derivative is zero.}
  \label{Fig17}
\end{figure}

\noindent
It follows therefore, as expected, that the infinite-order scaled kernel
function is differentiable. Note that, since the kernels on the right-hand
side have range $\epsilon/2$, and their points of application are distant
from each other by exactly $\epsilon$, each one is just outside the
support of the other. Therefore, the derivative is given by the
concatenation of two graphs just like the kernel itself, but with the
support scaled down from $\epsilon$ to $\epsilon/2$, with the amplitude
scaled up by the factor $2/\epsilon$, and with the sign of one of them
inverted. This is shown in Figure~\ref{Fig15}, containing a superposition
of the kernel and its rescaled first derivative.

\begin{figure}[ht]
  \centering
  \fbox{
    \epsfig{file=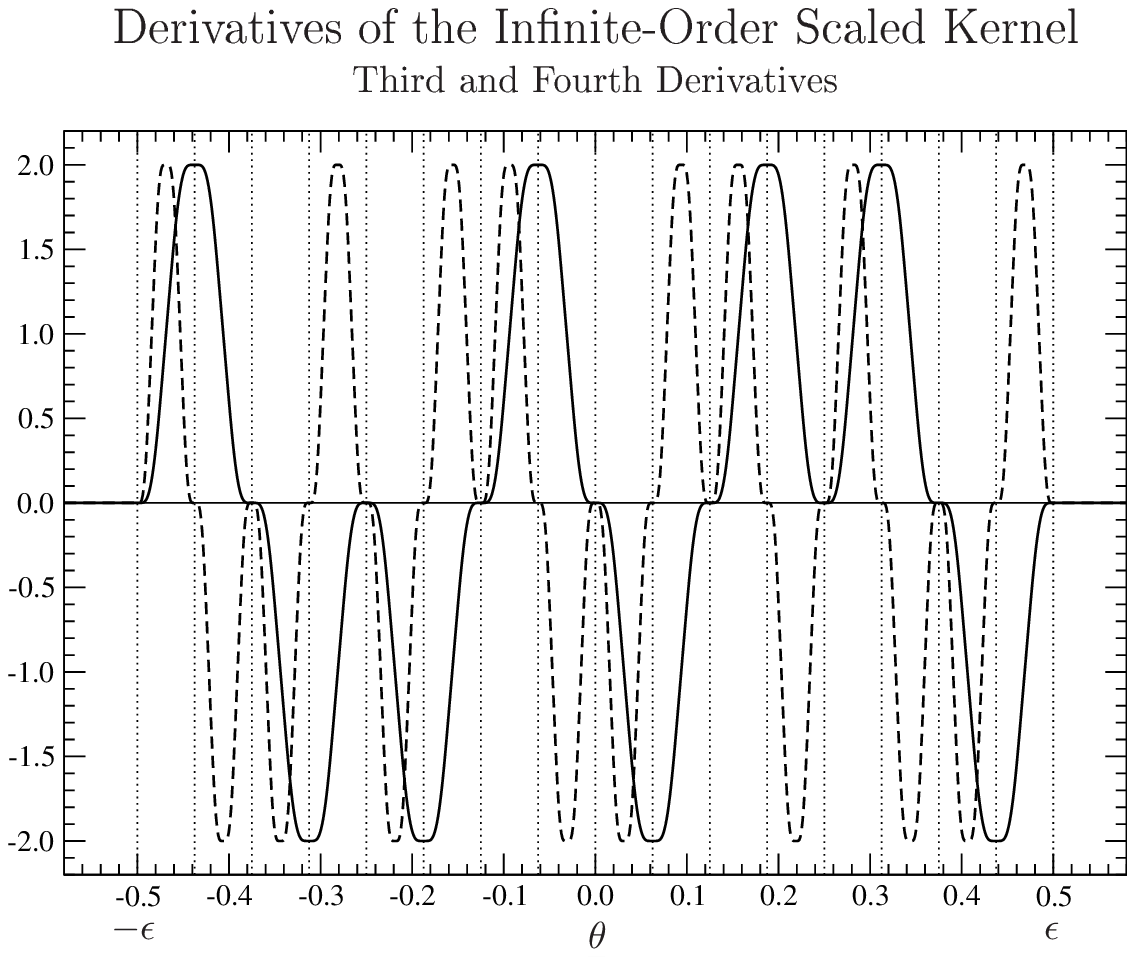,scale=1.0,angle=0}
  }
  \caption{The third (solid line) and fourth (dashed line) derivatives of
    the infinite-order scaled kernel, obtained via the use of their
    Fourier series, for $\epsilon=0.5$, for a large value of $N$ ($100$),
    plotted as functions of $\theta$ over the support interval
    $[-\epsilon,\epsilon]$. Both derivatives were rescaled down to have
    the same amplitude as the kernel. The dotted lines mark the points
    where the fourth derivative is zero.}
  \label{Fig18}
\end{figure}

One may now take one more derivative of the expression for the derivative
of the order-$N$ scaled kernel, thus obtaining an expression for the
corresponding second derivative. After that one may use again the relation
for the first derivative, thus iterating that relation, in order to obtain

\noindent
\begin{eqnarray*}
  \frac{d^{2}}{d\theta^{2}}
  \bar{K}_{\epsilon_{N}}^{(N)}(\theta)
  & = &
  \frac{1}{\epsilon}
  \left[
    \frac{d}{d\theta}
    \bar{K}_{\epsilon_{(N-1)}/2}^{(N-1)}\!\left(\theta+\epsilon/2\right)
    -
    \frac{d}{d\theta}
    \bar{K}_{\epsilon_{(N-1)}/2}^{(N-1)}\!\left(\theta-\epsilon/2\right)
  \right]
  \\
  & = &
  \frac{1}{\epsilon}
  \left\{
    \frac{2}{\epsilon}
    \left[
      \bar{K}_{\epsilon_{(N-2)}/4}^{(N-2)}\!\left(\theta+3\epsilon/4\right)
      -
      \bar{K}_{\epsilon_{(N-2)}/4}^{(N-2)}\!\left(\theta+\epsilon/4\right)
    \right]
  \right.
  +
  \\
  &   &
  \hspace{2em}
  -
  \left.
    \frac{2}{\epsilon}
    \left[
      \bar{K}_{\epsilon_{(N-2)}/4}^{(N-2)}\!\left(\theta-\epsilon/4\right)
      -
      \bar{K}_{\epsilon_{(N-2)}/4}^{(N-2)}\!\left(\theta-3\epsilon/4\right)
    \right]
  \right\}
  \\
  & = &
  \frac{2}{\epsilon^{2}}
  \left[
    \bar{K}_{\epsilon_{(N-2)}/4}^{(N-2)}\!\left(\theta+3\epsilon/4\right)
    -
    \bar{K}_{\epsilon_{(N-2)}/4}^{(N-2)}\!\left(\theta+\epsilon/4\right)
  \right.
  +
  \\
  &   &
  \hspace{2em}
  -
  \left.
    \bar{K}_{\epsilon_{(N-2)}/4}^{(N-2)}\!\left(\theta-\epsilon/4\right)
    +
    \bar{K}_{\epsilon_{(N-2)}/4}^{(N-2)}\!\left(\theta-3\epsilon/4\right)
  \right].
\end{eqnarray*}

\noindent
Once more the known existence of the $N\to\infty$ limit of the right-hand
side establishes the existence of the $N\to\infty$ limit of the left-hand
side, and therefore that the infinite-order scaled kernel is twice
differentiable. Taking the limit on both sides we get

\noindent
\begin{eqnarray*}
  \frac{d^{2}}{d\theta^{2}}
  \bar{K}_{\epsilon}^{\infty}(\theta)
  & = &
  \frac{2}{\epsilon^{2}}
  \left[
    \bar{K}_{\epsilon/4}^{\infty}(\theta+3\epsilon/4)
    -
    \bar{K}_{\epsilon/4}^{\infty}(\theta+\epsilon/4)
  \right.
  +
  \\
  &   &
  \hspace{2em}
  \left.
    -
    \bar{K}_{\epsilon/4}^{\infty}(\theta-\epsilon/4)
    +
    \bar{K}_{\epsilon/4}^{\infty}(\theta-3\epsilon/4)
  \right].
\end{eqnarray*}

\noindent
We now have four copies of the graph of the kernel, with range scaled down
to $\epsilon/4$ and amplitude scaled up by $2/\epsilon^{2}$, each one
outside the supports of the others, distributed in a regular way within
the interval of length $2\epsilon$ around $\theta$. This is shown in
Figure~\ref{Fig16}, containing a superposition of the rescaled first and
second derivatives of the kernel. As one can see in the subsequent
Figures~\ref{Fig17} and~\ref{Fig18}, the same type of relationship is also
true for all the higher-order derivatives. This is so because we can
iterate this relation indefinitely, so that any finite-order derivative of
$\bar{K}_{\epsilon}^{\infty}(\theta)$ can be written as a finite linear
combination of $\bar{K}_{\epsilon}^{\infty}(\theta)$ itself, with a
rescaled $\epsilon$ and a rescaled amplitude. Observe that this
constitutes independent proof that the infinite-order scaled kernel is a
$C^{\infty}$ function.

\subsection{Non-Analyticity of the Infinite-Order Scaled
  Kernel}\label{APPNonAnalytic}

From the construction described in the previous section for the order-$n$
derivatives of the infinite-order scaled kernel, which can all be written
in terms of the kernel itself, and from the fact that the infinite-order
scaled kernel is zero at the two ends of its support interval, it follows
at once that all the order-$n$ derivatives of the kernel are also zero at
these two points. Therefore the kernel and all its order-$n$ derivatives,
for all $n\geq 0$, are zero at the two ends $\theta=\pm\epsilon$ of the
support interval, as in fact we already knew, since this is also a
consequence of the fact that the kernel is a $C^{\infty}$ function over
its whole domain.

In a similar way, we may also determine other points within the support
interval where almost all the order-$n$ derivatives of the kernel are
zero. For example, at the central point, although the kernel itself is not
zero, we see that its first derivative is, and in fact iterating the
construction one can see that all the higher-order derivatives are zero
there. Therefore all the order-$n$ derivatives of the kernel, for all
$n\geq 1$, are zero at the central point $\theta=0$ of the support
interval. An examination of the situation at the two inflection points
$\theta=\pm\epsilon/2$ reveals that at those two points all the order-$n$
derivatives of the kernel, for all $n\geq 2$, are zero.

\begin{table}[ht]
  \centering
  \begin{tabular}{|c|c|c|}
    \hline
    Order  & Number of points & Null Derivatives \\
    \hline
    0      &              $2$ &        $m\geq 0$ \\
    1      &              $3$ &        $m\geq 1$ \\
    2      &              $5$ &        $m\geq 2$ \\
    3      &              $9$ &        $m\geq 3$ \\
    4      &             $17$ &        $m\geq 4$ \\
    \vdots &           \vdots &           \vdots \\
    n      &        $2^{n}+1$ &        $m\geq n$ \\
    \vdots &           \vdots &           \vdots \\
    \hline
  \end{tabular}
  \caption{Points where almost all derivatives of the infinite-order
    scaled kernel are zero. The integer $n$ is the order of the first null
    derivative, and $m$ gives the orders of all derivatives which are zero
    at the corresponding set of points.}
  \label{Tab01}
\end{table}

The iteration of this process of analysis can be continued indefinitely,
with the result that there are sets of increasing numbers of points
regularly spaced in the support interval where all derivatives above a
certain order are zero. This can be systematized as shown in
Table~\ref{Tab01}. We see therefore that there is a set of $2^{n}+1$
points regularly spaced within the support interval where all derivatives
with order $n$ or larger are zero. In the $n\to\infty$ limit this set of
points tends to be densely distributed within the support interval.
Outside the support interval the derivatives of all orders are zero at all
points, of course, since the kernel is identically zero there.

In we assemble the Taylor series of the kernel function around one of the
points where all the derivatives of order $n$ and larger are zero, we
obtain a convergent power series, which is in fact a polynomial of order
$n-1$. Since the kernel function is obviously not such a polynomial, it is
therefore {\em not} represented by its convergent Taylor series around
this reference point, at any points other than the reference point itself.
Since in order to be analytic the kernel function would have to be so
represented within an open set around the reference point, it follows that
it is not analytic at any of these points. Since this set of points tends
to become densely distributed within the support interval, we may conclude
that the kernel function is {\em not} analytic at all points of the
support interval.

One can try to extend this argument to show in a somewhat heuristic way
that the kernel function cannot be represented by a convergent power
series around any point of the support interval, whether or not it is in
the dense subset. Let us consider a point where the kernel function has
non-zero derivatives of arbitrarily high orders, and where the Taylor
series built from them converges in an open neighborhood of that reference
point. This implies that the point at issue is not in the dense subset.
Note that since the kernel function is $C^{\infty}$ we know that all its
derivatives at the point exist, whether or not they are zero. Since the
subset of points discussed above is dense in the support interval, there
is at least one point of the dense subset within this open neighborhood.
Therefore there is another Taylor series around this point, which is also
convergent.

Since both are Taylor series of the same function and converge in a common
domain, we must be able to transform each one into the other by a
transformation of coordinates that is a simple shift of the argument of
the series. Note that all the derivatives of the kernel function, of all
orders, are themselves continuous and differentiable functions. However,
no such transformation of variables can transform the second series, which
is a polynomial of finite order, into a series such as the first one, with
no upper bound to the powers present in it. This seems to produce an
absurd situation. Therefore, one is led to think that either the first
series cannot be a convergent series, or it must converge to some function
other than the kernel function. In any case, it follows that the kernel
function is not represented by this Taylor series either, and once again
that it cannot be analytic at the point under discussion.

One may wonder about whether the real infinite-order kernel function can
be extended analytically to the complex plane. It is clear that this
cannot be done in the usual way, with the simple exchange of its argument
by a complex variable. In addition to this, we know that it can be
obtained as the limit to the unit circle of an inner analytic function,
and that the inner analytic function has a densely distributed set of
singularities on the support of the kernel. It is therefore reasonable to
think that this is not possible, but no complete proof of this is
currently available.

\end{document}